\documentclass[11pt, reqno]{amsart}

\usepackage[mathscr]{eucal}
\usepackage{mathrsfs}
\usepackage{xypic}
\usepackage{graphicx}
\usepackage{floatrow}
\usepackage[font=scriptsize]{caption}
\usepackage{amsfonts}
\usepackage{amsmath}
\usepackage{amsthm}
\usepackage{amssymb}
\usepackage{latexsym}
\usepackage{pict2e}
\usepackage[all]{xy}
\usepackage{hyperref}
\usepackage{baskervald}

\allowdisplaybreaks

\newtheorem{theorem}{Theorem}
\newtheorem{lemma}[theorem]{Lemma}

\theoremstyle{definition}

\newtheorem{example}[theorem]{Example}

\newtheorem{definition}[theorem]{Definition}
\newtheorem{definition-lemma}[theorem]{Definition-Lemma}
\newtheorem{example-lemma}[theorem]{Example-Lemma}
\newtheorem{assumption}[theorem]{Assumption}

\newtheorem*{convention}{Convention}
\newtheorem*{ack}{Acknowledgements}

\newcommand{\HH}{{\mathrm{H}}}

\newcommand{\R}{\mathbb{R}}
\newcommand{\Z}{{\mathbb{Z}}}

\newcommand{\A}{\mathcal{A}}

\newcommand{\ov}{\overline}

\newcommand{\Hom}{\mathrm{Hom}}
\newcommand{\Hoch}{\mathrm{Hoch}}
\newcommand{\Cycl}{\mathrm{Cycl}}
\newcommand{\Cont}{\mathrm{Cont}}
\newcommand{\Fuk}{\mathcal Fuk}
\newcommand{\Lag}{{Lag}}

\begin{document}

\title[Cyclic homology of Fukaya categories and the linearized contact homology]
{\Large Cyclic homology of Fukaya categories and\\ the linearized contact homology}

\author{
Xiaojun Chen,  
Hai-Long Her, 
Shanzhong Sun 
}

\thanks{
H.-L. Her was partially supported by NSFC (No.10901084); 
S. Sun was partially supported by NSFC (No.10731080, 11131004), 
PHR201106118, the Institute of Mathematics and Interdisciplinary Science at CNU
}

\address{X. Chen: School of Mathematics, Sichuan University, Chengdu 610064 P. R. China}
\email{xjchen@scu.edu.cn}

\address{H.-L. Her: School of Mathematical Sciences, Nanjing Normal University, Nanjing 210046 P. R. China}
\email{hailongher@126.com}

\address{S. Sun: School of Mathematics, Capital Normal University, Beijing 100048 P. R. China}
\email{sunsz@mail.cnu.edu.cn}

\begin{abstract}
Let $M$ be an exact symplectic manifold with contact type boundary such that $c_1(M)=0$. 
In this paper we show that
the cyclic cohomology of the Fukaya category of $M$ has the structure of an involutive Lie bialgebra.
Inspired by a work of Cieliebak-Latschev we show that there is a Lie bialgebra homomorphism
from the linearized contact homology of $M$ to the cyclic cohomology of the Fukaya category.
Our study is also motivated by string topology and 2-dimensional topological conformal
field theory.
\end{abstract}

\maketitle

\setcounter{tocdepth}{1}

\tableofcontents

\section{Introduction}

\renewcommand{\thetheorem}{\Alph{theorem}}

During the past two decades,
great achievements have been obtained
in the understanding of the geometry and topology
of symplectic manifolds. Among them are the Fukaya
category of symplectic manifolds, originated by Fukaya (\cite{Fukaya93}),
and the symplectic field theory (SFT), introduced by Eliashberg-Givental-Hofer (\cite{EGH}).
In this paper,
we study the cyclic homology of the Fukaya category
of an exact symplectic manifold with contact type boundary
and its relations to the linearized contact homology
defined in SFT.

The main motivation of our study is
Kontsevich and Costello's work on topological conformal field theories
(\cite{Kontsevich1,Kontsevich} and \cite{Costello2}),
and the work
of Cieliebak-Latschev (\cite{CL}), which studies
the algebraic structures on the linearized contact homology and their relations with string topology.
Let us explain in more detail.

\subsection{}
In \cite{Kontsevich} Kontsevich associates to an A$_\infty$ algebra with
a cyclically invariant inner product a cohomology class on
the compactified moduli spaces of marked Riemann surfaces, and from this as well
as some analogous examples he first raised his theory of noncommutative symplectic geometries.
Such an A$_\infty$ algebra was later called by him a {\it Calabi-Yau A$_\infty$ algebra},
and was systematically studied in his work with Soibelman (\cite{KS}).
In his talk at the Hodge Centennial Conference (\cite{Kontsevich}) he showed that a Calabi-Yau
A$_\infty$ algebra, or more generally, a Calabi-Yau A$_\infty$ algebra ``with several objects'', {\it i.e.} a
Calabi-Yau A$_\infty$ category, is an {\it open} topological
conformal field theory (TCFT), and its Hochschild cohomology
is in fact a {\it closed} topological conformal field theory.
As a particular example, he conjectured that the Fukaya category of a symplectic manifold is a Calabi-Yau
category.
Similar restults have been
obtained by Costello in \cite{Costello2}.
We refer the reader to Getzler \cite{Getzler,Getzler2} and Costello \cite{Costello2}
for more details and interesting results about TCFT's. 
Inspired by \cite{Kontsevich} and string topology (\cite{CS02}), the first
author showed that the cyclic cohomology of a Calabi-Yau A$_\infty$ category endows the
structure of an involutive Lie bialgebra (\cite{Chen}).

While the existence of the cyclically invariant inner product on the Fukaya category is still
under verification (for some partial results see \cite{Fukaya2}),
the work of Kontsevich and Costello remains a guiding philosophy for our study.
Our first theorem in the following says that the Lie bialgebra exists on the cyclic cohomology
of the Fukaya category (no non-degenerate pairing is assumed):

\begin{theorem}[Theorem \ref{Liebi_Fuk}]\label{theoremA}
Let $M$ be an exact symplectic manifold with contact type boundary such that $c_1(M)=0$.
Then the cyclic cohomology of Fukaya category $\mathcal Fuk(M)$ of $M$ has the structure of an
involutive Lie bialgebra.
\end{theorem}

The key point in the above theorem is that, in the construction of the Fukaya category,
counting the pseudo-holomorphic disks is {\it a priori} cyclically invariant.

Let $M$ be as in Theorem \ref{theoremA}. Denote its boundary by $W$.
The symplectic field theory of Eliashberg-Givental-Hofer
relates the geometry of closed Reeb orbits on $W$ with the geometry of pseudo-holomorphic
curves on the symplectic completion $\widehat M$ of $M$.
In fact, pseudo-holomorphic curves in $\widehat M$ asymptotic to closed Reeb orbits
share a lot of properties of a closed TCFT; for more details, see \cite{EGH}. Among all
the interesting properties of SFT,
Cieliebak-Latschev (\cite[Theorem A]{CL}) proved that the
{\it linearized contact homology}  of $M$, denoted by $\mathrm{CH}^{\mathrm{lin}}_*(M)$,
which basically arises from counting
the pseudo-holomorphic cyclinders in $\widehat M$, is in fact
a Lie bialgebra. With the guiding philosophy of Kontsevich and Costello that
the cyclic cohomology of the Fukaya category of $M$ is a closed TCFT (and in fact,
this closed TCFT is universal in the sense that any other closed TCFT will factor through it), one
might wonder whether there is a direct relationship between them.
The following theorem gives an affirmative answer to this question, which
is also inspired by a theorem of Cieliebak-Latschev
in the same article (\cite[Theorem B]{CL}) and by Seidel \cite{Seidel02}, and may be viewed as a
realization of the Holographic Principle (or in some other words, the Bulk-Boundary Correspondence) in physics:

\begin{theorem}[Theorem \ref{thm_liebimap_comp}]\label{mapfromcont}
Let $M$ be an exact symplectic manifold with contact type boundary such that $c_1(M)=0$.
There is a chain map from the linearized contact complex $\Cont_*^{\mathrm{lin}}(M)$
to the cyclic cochain complex $\Cycl^*(\Fuk(M))$ of $\Fuk(M)$:
\begin{equation}\label{ContFuk}
f: \Cont^{\mathrm{lin}}_*(M)\longrightarrow \Cycl^*(\Fuk(M)),
\end{equation}
which induces a Lie bialgebra homomorphism on their homology.
\end{theorem}

The homomorphism in the above theorem is similar to the one given in \cite{CL,Seidel02}, which
is given by counting the pseudo-holomorphic cylinders with one end approaching a closed Reeb orbit and the other end
lying in several Lagrangian submanifolds. We remark that  the study of pseudo-holomorphic
curves with some boundary components approaching the closed
Reeb orbits and some other lying in Lagrangian submanifolds
(each boundary component lying exactly in {\it one} Lagrangian submanifold) has already been discussed in \cite{EGH}.
In the case of pseudo-holomorphic cylinders, Cieliebak-Latschev first showed that
they in fact induce
a homomorphism of Lie bialgebras from the linearized contact homology to the $S^1$-equivariant homology of the
free loop space of the Lagrangian submanifold being considered.
However, it is hard to see that the homomorphism they discovered is {\it on the chain level}.

We remark that, in fact, on the chain level the homomorphism in Theorem~\ref{mapfromcont}
is a {\it Lie$_\infty$ homomorphism}
(or more generally, a {\it BV$_\infty$ homomorphism} in the sense of Cieliebak-Latschev). More details can be
found in \S\ref{sect_rel}.
The key point of Theorem~\ref{mapfromcont} is that, when considering several Lagrangian submanifolds together,
{\it i.e.} the Fukaya category, we do get a homomorphism on the chain level from the linearized
contact chain complex to the cyclic cochain complex of the Fukaya category.
This idea seems to have first appeared in Seidel \cite{Seidel02}, where he
considered a version of pseudo-holomorphic cylinders,
with one boundary component being a Hamiltonian closed orbit and the other 
lying on several Lagrangian submanifolds.

Seidel's homomorphism is from the symplectic homology to the Hochschild cohomology
of the Fukaya category. Our construction may be viewed as a cyclic version of his.
On the other hand, from the work of Bourgeois-Oancea (\cite{BO}) there is a deep relationship between
the symplectic homology and the linearized contact homology; in fact, there is a long exact
sequence which relates them. As is probably well-known
from cyclic homology theory and non-commutative geometry,
there is also a long exact sequence, called the Connes exact sequence, connecting
the Hochschild homology and cyclic homology, too.
In fact, Jones showed in \cite{Jones}
that for a simply connected manifold,
the Hochschild homology (resp. the cyclic homology)
of its de Rham complex is isomorphic to the cohomology (resp. the $S^1$-equivariant cohomology)
of its free loop space.
Such a theorem is also implicit (or almost explicit, as we would say) in the work of K.-T. Chen \cite{KTChen},
and a good reference for it is Getzler-Jones-Petrack \cite{GJP}. 

All these results in fact fit into a package, called {\it string topology}, initiated by Chas-Sullivan (\cite{CS99}).
Roughly speaking, string topology is a study of the topological structures on the free loop space
of compact manifolds. 
Chas-Sullivan proved in \cite{CS02} that
the $S^1$-equivariant homology of the free loop space of a compact manifold relative to the
constant loops has the structure of a Lie
bialgebra. Such a Lie bialgebra, as shown by Cieliebak-Latschev (\cite[Theorem C]{CL}),
is isomorphic to the linearized contact homology of the cotangent bundle of the manifold.

\subsection{Some related works}

Besides the articles that we have cited above, there are several other works
that are related to the interest of this paper.
First, we have benefited a lot from the paper of Abbondandolo-Schwarz \cite{AS},
which gives a complete treatment of the isomorphism between
the Floer/symplectic homology of cotangent bundles and the loop product defined in
string topology by Chas-Sullivan in \cite{CS99}.
We are also inspired by the work of Bourgeois-Oancea \cite{BO}.
As we have said above, they have shown that there is a long exact sequence
$$
\xymatrix{
\cdots\ar[r]&\mathrm{SH}_*(M)\ar[r]&\mathrm{CH}_*^{\mathrm{lin}}(M)\ar[r]&
\mathrm{CH}_{*-2}^{\mathrm{lin}}(M)\ar[r]&
\mathrm{SH}_{*-1}(M)\ar[r]&\cdots
}
$$
where $\mathrm{SH}_*(-)$ is the symplectic homology and $\mathrm{CH}^{\mathrm{lin}}_*(-)$
is the linearized contact homology, respectively.
From the work of Seidel \cite{Seidel02}
and our above theorem,
there should be the following morphism of long exact sequences
$$
\xymatrixcolsep{1.5pc}
\xymatrix{
\ar[r]&\mathrm{SH}_*(M)\ar[r]\ar[d]&\mathrm{CH}_*^{\mathrm{lin}}(M)\ar[r]\ar[d]&
\mathrm{CH}_{*-2}^{\mathrm{lin}}(M)\ar[r]\ar[d]&
\mathrm{SH}_{*-1}(M)\ar[r]\ar[d]&\cdots\\
\ar[r]&\mathrm{HH}^*(\Fuk(M))\ar[r]&
\mathrm{HC}^*(\Fuk(M))\ar[r]&
\mathrm{HC}^{*-2}(\Fuk(M))\ar[r]&
\mathrm{HH}^{*-1}(\Fuk(M))\ar[r]&\cdots
}$$
where $\mathrm{HH}^*(-)$ is the Hochschild cohomology with values in a ground field of
characteristic zero and
$\mathrm{HC}^*(-)$ is the cyclic cohomology, respectively,
and the long exact in the bottom line
is the Connes exact sequence for A$_\infty$ categories.



\indexspace

So far we have mentioned in this paper several homology theories:
symplectic homology, linearized contact homology, Hochschild cohomology,
cyclic cohomology, $\cdots$, and various homomorphisms among them.
These homomorphisms are usually given by counting the pseudo-holomorphic
cylinders in each situation being considered.
We list the types of pseudo-holomorphic cylinders that
are studied by different authors in literature, which may
help the reader
to get some more idea on these constructions.

\begin{figure}[h]\label{figures}
\begin{floatrow}[2]
\ffigbox[8.6cm]{\caption{Homomorphism
from symplectic homology to the Hochschild cohomology of
Fukaya category (Seidel \cite{Seidel02}).}}
{\includegraphics[height=3.6cm]{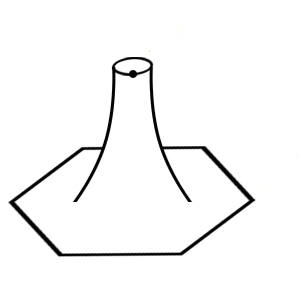}}
\ffigbox[8.6cm]{\caption{Homomorphism
from linearized contact homology to the $S^1$-equivariant
homology of the free loop space of the zero section (Cieliebak-Latschev \cite{CL}).}}
{\includegraphics[height=3.6cm]{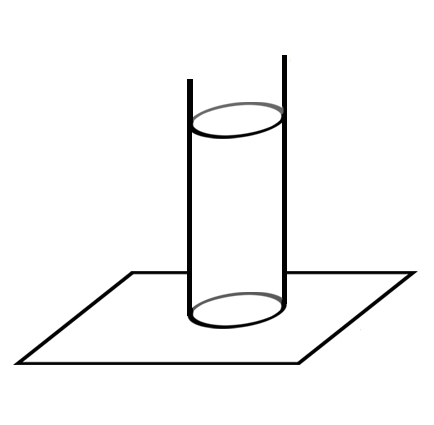}}
\end{floatrow}

\begin{floatrow}[2]
\ffigbox[8.6cm]{\caption{Homomorphism
from $S^1$-parametrized
linearized contact homology to
symplectic homology (Bourgeois-Oancea \cite{BO}).}}
{\includegraphics[height=3.6cm]{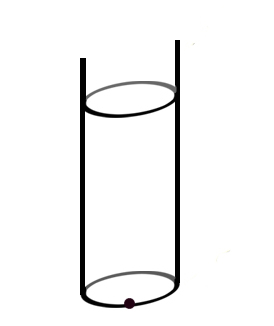}}
\ffigbox[8.6cm]{\caption{Homomorphism
in current paper from
linearized contact homology to the cyclic cohomology of Fukaya category.}}
{\includegraphics[height=3.6cm]{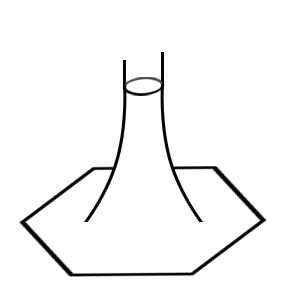}}
\end{floatrow}
\end{figure}

Also, after the first draft was posted on arXiv, we are informed by Eliashberg that 
he, Bourgeois and Ekholm have obtained some results analogous to ours (\cite{BEE,BEE2}).
In these two papers they study a variety of algebraic structures (including Lie bialgebra) on the Legendrian homology,
and relate them to symplectic homology (the pictures are similar to above). We are grateful
to him for acknowledging us about their work.

\subsection{}
The rest of the paper is devoted to the proof of the above theorems. It is organized as follows:
In \S2, we first collect some facts on the Fukaya category of exact symplectic manifolds, and then prove
Theorem \ref{theoremA} in \S\ref{sect_liebialg}.
In \S\ref{sect_lin} we recall the definition and some properties of linearized contact homology.
In \S\ref{sect_ana} we set up the necessary analytic machinery for the moduli space of pseudo-holomorphic
curves to be considered later.
In \S\ref{sect_rel} we prove Theorem \ref{mapfromcont}.
In \S\ref{sect_cot}, we study the example of cotangent bundles
and relate it with string topology.

\subsection{}
Finally, we remark that in this article,
TCFT and string topology have served as an inspiring motivation of our study.
However, the main body of this article is independent of these two theories, and we shall not discuss 
any details of them
in the rest of the paper. The interested reader may refer to the literature.

\begin{ack}
Some results in this paper have been presented by the first author at the Chern Institute, Tianjin in 2010
and on the Sullivanfest, Stony Brook in 2011; he would like to thank these two places for their
hospitalities. The first author also thanks
Nanjing Normal University, the Capital Normal University
and NCTS in Taiwan for their supports and hospitalities, and Liang Kong for very inspiring conversations.
The second author would like to thank BICMR and the Capital Normal University for their hospitalities, 
and Gang Tian for helpful suggestions.
All three authors thank Yiming Long and Yongbin Ruan for their encouragements during
these years, and D. Pomerleano for pointing out an error in previous draft. 
\end{ack}

\setcounter{theorem}{0}
\renewcommand{\thetheorem}{\arabic{theorem}}

\section{Fukaya category of exact symplectic manifolds}

\subsection{A$_\infty$ categories}

\begin{definition}[A$_\infty$ Category]\label{def_Ainfty}
An {\it A$_\infty$ category} $\mathcal A$ consists of a set of objects $\mathcal Ob(\mathcal A)$,
a graded vector space
$\Hom_\A(A_1,A_2)$ for each pair of objects $A_1,A_2\in\mathcal Ob(\mathcal A)$, and a sequence of
operators:
$$m_n^\A: \Hom_\A(A_1,A_2)\otimes \Hom_\A(A_2, A_3)
\otimes\cdots\otimes \Hom_\A(A_{n},A_{n+1})\to \Hom_\A(A_1,A_{n+1}), $$
where $|m_n^\A|=2-n$, for $n=1,2,\cdots,$ satisfying the following A$_\infty$ relations:
\begin{equation}\label{higher_htpy}
\sum_{p=1}^{n}\sum_{k=1}^{n-p+1}(-1)^{\mu_{npk}} m_{n-k+1}^\A(a_1,\cdots,a_{p-1},
m_{k}^\A(a_{p},\cdots,a_{p+k-1}),a_{p+k},\cdots,a_n)=0,
\end{equation}
where $a_i\in\Hom_\A(A_{i},A_{i+1})$, for $i=1,2,\cdots,n$, and
$\mu_{npk}=\sum\limits_{r=1}^n(r-1)|a_r|+\sum\limits_{s=1}^{p-1} k |a_s|+p(2-k)$.
\end{definition}

If an A$_\infty$ category has one object, say $A$, then $\Hom_\A(A,A)$ is an
{\it A$_\infty$ algebra}; and if furthermore, all $m_i$, $i\ge 3$, vanish,
then $\A$ is the usual differential graded (DG) algebra.

\begin{convention}[The Sign Rule]
The sign in equation (\ref{higher_htpy}) is always complicated. The rule is given as follows.
First, for a graded vector space $V$, let $\ov V$ be
the desuspension of $V$, that is, $(\ov V)_i=V_{i+1}$. Let $\Sigma: V\to \ov V$
be the identity map which maps $v$ to $\ov v$, and let
$$
\begin{array}{cccl}
\Sigma^{\otimes n}:&V\otimes\cdots\otimes V& \longrightarrow&\ov V\otimes\cdots\otimes\ov V\\
&v_1\otimes\cdots\otimes v_n&\longmapsto&(-1)^{(n-1)|v_1|+\cdots+|v_{n-1}|}\ov v_1\otimes\cdots\otimes \ov v_n
\end{array}
$$
be the $n$-folder tensor of $\Sigma$.
Let $\ov m_n:(\ov V)^{\otimes n}\to \ov V$ be the degree 1 map such that the diagram
\begin{equation}\label{signrule}
\xymatrixcolsep{5pc}
\xymatrix{V\otimes \cdots\otimes V\ar[d]^{\Sigma^{\otimes n}}
\ar[r]^-{m_n}&V\ar[d]^\Sigma\\
\ov V\otimes\cdots\otimes\ov V\ar[r]^-{\ov m_n}&\ov V
}\end{equation}
commutes. Then equation (\ref{higher_htpy}) is nothing but
\begin{equation}\label{htpyrel}
\sum_{p=1}^{n}\sum_{k=1}^{n-p+1}(-1)^{|\ov a_1|+\cdots+|\ov a_{p-1}|}\ov m_{n-k+1}
(\ov a_1,\cdots,\ov a_{p-1}, \ov m_k(\ov a_{p},\cdots, \ov a_{p+k-1}),
\ov a_{p+k},\cdots,\ov a_n)=0.
\end{equation}
The sign that appears in equation (\ref{htpyrel}) follows from the usual Koszul convention rule.
Namely, the canonical isomorphism
$V\otimes W\stackrel{\cong}{\to}W\otimes V$
is given by
$a\otimes b\mapsto (-1)^{|a||b|}b\otimes a$.
One then obtains equation (\ref{higher_htpy}) by converting equation (\ref{htpyrel}) via diagram (\ref{signrule}).
In the following all signs are assigned in this way,
and therefore we will just simply write $\pm$, without specifying their particular value.
\end{convention}

In the following, if $\A$ is clear from the context, we will simply write
$\Hom_\A(A,B)$ and $m_i^\A$ as $\Hom(A,B)$ and $m_i$.
For an A$_\infty$ category $\A$, since $m_1\circ m_1=0$,
one obtains the {\it cohomological category} $\HH(\A)$ of $\A$. Namely,
the objects of $\HH(\A)$ are the same as the objects of $\A$ while the morphisms from $A$ to $B$ are the
cohomology classes $\HH(\Hom(A,B),m_1)$. $\HH(\A)$ is a not-necessarily-unital category.

\begin{definition}[Hochschild Homology]\label{def_hoch}
Let $\mathcal A$ be an A$_\infty$ category. The {\it Hochschild
chain complex} $\Hoch_*(\mathcal A)$ of $\mathcal A$ is the chain
complex whose underlying vector space is
\begin{equation}\label{defHoch}
\bigoplus_{n=1}^{\infty}\bigoplus_{A_1,A_2,\cdots,A_{n+1}\in\mathcal Ob(\mathcal A)}
\ov{\Hom}(A_1,A_2)\otimes\ov{\Hom}(A_{2}, A_3)\otimes\cdots\otimes\ov{\Hom}(A_{n+1},A_1)
\end{equation}
with differential $b$ defined by
\begin{eqnarray*}&&b(\ov a_1,\ov a_2,\cdots,\ov a_{n+1})\\
&=&\sum_{k=1}^{n+1}\sum_{i=1}^{n-k+1}\pm(\ov a_1,\cdots, \ov a_{k-1}, \ov m_i(\ov a_k,
\cdots,\ov a_{k+i-1}),\ov a_{k+i},\cdots, \ov a_{n+1})\\
&+&\sum_{j=0}^{n-1}\sum_{i=1}^{n-j}
\pm(\ov m_{i+j+1}(\ov a_{n-j+1}, \ov a_{n-j+2},\cdots,\ov a_{n+1}, \ov a_{1},\ov a_{2},\cdots,\ov a_{i}),
\ov a_{i+1}\cdots,\ov a_{n-j}).
\end{eqnarray*}
The associated homology is called the {\it Hochschild homology} of $\mathcal A$,
and is denoted by $\mathrm{HH}_*(\mathcal A)$.
\end{definition}

\begin{definition}[Cyclic Homology]\label{def_cycl}
Suppose $\mathcal A$ is an A$_\infty$ category.
Let
\begin{multline*}t_{n}: \ov{\Hom}(A_1,A_2)\otimes\ov{\Hom}
(A_2, A_3)\otimes\cdots\otimes\ov{\Hom}(A_{n+1},A_1)\\
\longrightarrow
\ov{\Hom}(A_{n+1},A_1)\otimes\ov{\Hom}(A_1,A_2)\cdots\otimes\ov{\Hom}(A_{n},A_{n+1}),
\end{multline*}
for $n=0,1,2, \cdots$, be the linear map
\begin{equation}\label{def_t}
t_n(\ov a_1, \ov a_{2},\cdots,\ov a_{n+1}):=
(-1)^{|\ov a_{n+1}|(|\ov a_1|+\cdots+|\ov a_{n}|)}(\ov a_{n+1},\ov a_1,\cdots, \ov a_{n}).
\end{equation}
Extend $t_n$ to $\Hoch_*(\mathcal A)$ trivially, and let
$t=t_0+t_1+t_2+\cdots$, the cokernel
$\Hoch_*(\mathcal A)/(1-t)$ of $1-t$ forms a chain complex with the induced differential from the
Hochschild complex (still denoted by $b$).
Such chain complex is denoted by $\Cycl_*(\mathcal A)$, and is called the {\it Connes cyclic
complex} of $\mathcal A$. Its homology is called the {\it cyclic homology} of $\mathcal A$, and is denoted by
$\mathrm{HC}_*(\mathcal A)$.

The {\it cyclic cohomology} of $\mathcal A$ is the homology of the dual cochain complex $\Cycl^*(\mathcal A)$
of $\Cycl_*(\mathcal A)$. Namely, suppose
$f\in\Hom(\Hoch_*(\mathcal A), k)$, then $f\in\Cycl^*(\mathcal A)$ if and only if for all
$\alpha\in\Hoch_*(\mathcal A)$, $f(\alpha)=f(t(\alpha))$.
\end{definition}

\subsection{Review of the Fukaya category}\label{fukcat}

In this subsection we briefly recall the construction of the Fukaya category in
{\it exact} symplectic manifolds. We adopt the setting of Seidel \cite{Seidel}.
All details and proofs are omitted. The construction
in a general symplectic manifold is given in Fukaya \cite[Chapter 1]{Fukaya}, largely based on
the work \cite{FOOO}; however,
we do not need to be such general.

An {\it exact symplectic manifold with contact type boundary} is a
quadruple $(M,\omega,\eta,J)$, where $M$ is a compact $2n$ dimensional manifold with
boundary, $\omega$ is a symplectic 2-form on $M$, $\eta$ is a
1-form such that $\omega=d\eta$ and $J$ is a $\omega$-compatible
almost complex structure. These data also satisfy the
following two convexity conditions:
\begin{itemize}\item The negative Liouville
vector field defined by
$\omega(\cdot,X_\eta)=\eta$
points strictly inwards along the boundary of $M$;
\item The boundary of $M$ is weakly $J$-convex, which means that any
pseudo-holomorphic curves cannot touch the boundary unless they are
completely contained in it.
\end{itemize}
An $n$ dimensional submanifold $L\subset M$ is
called {\it Lagrangian} if $\omega|L=0$. We always assume $L$ is
{\it closed} and is disjoint from the boundary of $M$. $L$ is called {\it
exact} if $\eta|L$ is an exact 1-form.

\begin{assumption}\label{assum-Fuk}
In the following, for a symplectic manifold $M$ with or without (contact type) boundary, we shall always assume $c_1(M)=0$,
and for a Lagrangian submanifold $L$ in $M$, we shall always assume it is {\it admissible}, namely,
(1) $\eta|_L$ is exact;
(2) $L$ has vanishing Maslov class; and
(3) $L$ is spin.
\end{assumption}

\begin{example}[Cotangent Bundles]
Let $N$ be a simply connected, compact spin manifold. Let $T^*N$ be the cotangent bundle of $N$
with the canonical symplectic structure. The cotangent disk bundle of $N$ is an exact symplectic
with contact type boundary. In particular, $N$, viewed as the zero section of $T^*N$, is an admissible Lagrangian submanifold.
\end{example}

Intuitively, the Fukaya category $\mathcal Fuk(M)$ of $M$
is defined as follows: the objects are the admissible Lagrangian submanifolds;
suppose $L_1,L_2$ are two objects,
$\Hom(L_1,L_2)$, called the {\it Floer cochain complex},
is spanned by the transversal intersection points of $L_1$ and $L_2$, and
for $n$ objects $L_1,\cdots, L_{n+1}$,
$$\ov m_n:\ov{\Hom}(L_1,L_2)\otimes\ov{\Hom}
(L_{2},L_3)\otimes\cdots\otimes\ov{\Hom}(L_{n},L_{n+1})\to\ov{\Hom}(L_1,L_{n+1})$$
is given by counting pseudo-holomorphic disks
whose boundary lying in $L_1,L_2,\cdots, L_{n+1}$.
More precisely, if $a_1\in\Hom(L_1,L_2),\cdots,a_n\in\Hom(L_n,L_{n+1})$,
$$\ov m_n(\ov a_1,\ov a_2,\cdots,\ov a_n)=
\sum_{a\in L_1\cap L_{n+1}}\#(\mathcal M(a,a_1,\cdots,a_n))\cdot\ov a,$$
where $\mathcal M(a,a_1,\cdots,a_n)$ is the moduli space of pseudo-holomorphic disks
with $n+1$ (anti-clockwise) cyclically ordered marked points in its boundary, such that
these marked points are mapped onto $a,a_{n},\cdots,a_1$ and that the
rest of the boundary lie in $L_1,L_2,\cdots, L_{n+1}$.
The A$_\infty$ relations (equation (\ref{higher_htpy}))
follow from the compactification of $\mathcal M(a,a_1,\cdots,a_n)$
where those pseudo-holomorphic disks with all possible
``bubbling-off'' disks are added.

This is a very rough description of the construction of the Fukaya category.
It is only partially defined in the sense
that we have assumed that all Lagrangian submanifolds are transversal;
also, the Floer cochains thus described are only $\mathbb Z_2$ graded.
To make the Fukaya category be fully defined
and be graded over $\mathbb Z$, we have to introduce the following concepts. Let us
do it one by one.

\subsubsection{Pointed-boundary Riemann surfaces}\label{PbRS}

Suppose $\widehat S$ be a Riemann surface with boundary, and $\Sigma$ is
a set of boundary points in $\widehat S$.
$\Sigma$ is divided into two subsets $\Sigma^+$ and $\Sigma^{-}$,
called the output subset and input subset. Now to each $\zeta\in\Sigma$, one associates:
\begin{itemize}
\item
two admissible Lagrangian submanifolds $(L_{\zeta 0},L_{\zeta 1})$, where $L_{\zeta 0}$
is uniquely attached to the boundary component of $S:=\widehat S\backslash\Sigma$ that
comes before $\zeta$ and $L_{\zeta 1}$ is uniquely attached to the boundary component
(with induced orientation) that comes
after $\zeta$ if $\zeta\in\Sigma^+$; otherwise if $\zeta\in\Sigma^{-}$,
$L_{\zeta 1}$ comes before $L_{\zeta 0}$. These Lagrangian submanifolds
are called the {\it Lagrangian labels}.
\item   

a  {\it strip-like end} which is a proper holomorphic embedding $\epsilon_\zeta:Z^\pm\to S$ satisfying
\begin{eqnarray}
 \epsilon^{-1}_\zeta(\partial S)=\R^\pm\times\{0,1\}\ \ \ &{\rm and}&\ \ \ \lim_{s\to\pm\infty}\epsilon_\zeta(s,\cdot)=\zeta,
\end{eqnarray}
where $Z^\pm=\R^\pm\times [0,1]$ denotes the semi-infinite strips. If $\Sigma$ consists of more than one point, we also need
the additional requirement that the images of the $\epsilon_\zeta$ are pairwise disjoint.
Such an $S$ is called a {\it pointed-boundary Riemann surface with strip-like ends}.
\end{itemize}

\subsubsection{Floer data and perturbation data}

Suppose $(M,\omega, J)$ is an exact symplectic manifold with contact type boundary.
Let $\mathscr J$ be the space of all $\omega$-compatible almost complex structures on $M$
which agree with the given $J$ near the boundary; and let $\mathscr H=C_c^\infty(int(M),\mathbb R)$
be the space of smooth functions on $M$ vanishing near the boundary.

\begin{definition}[Floer Datum]\label{F-data}
For each {\it ordered} pair of Lagrangian submanifolds
$L_0, L_1\subset M$, a {\it Floer datum} consists of
$H_{L_0, L_1}\in C^\infty([0,1], \mathscr H)$ and $J\in C^\infty([0,1],\mathscr J)$, with the
following property: if $X$ is the time-dependent Hamiltonian vector field of $H$ and
$\phi$ is its flow, then $\phi^1(L_0)$ intersects $L_1$ transversally.
\end{definition}

\begin{definition}[Perturbation Datum]
Let $S$ be a pointed-boundary Riemann surface with Lagrangian labels. Suppose we have
chosen strip-like ends for it, and also a Floer datum $(H_\zeta, J_\zeta)$ for each
of the pairs of sub manifolds $(L_{\zeta 0}, L_{\zeta 1})$ associated to the points at infinity
$\zeta\in\Sigma$. A {\it perturbation datum} for $S$ is a pair $(K, J)$ where
\begin{itemize}\item $K\in\Omega^1(S,\mathscr H)$ satisfies
$K(\xi)|L_C=0$, for all $\xi\in TC$, where $C$ is a component of $\partial S$, and

\item $J$ is a family of almost complex structures $J\in C^\infty(S,\mathscr J)$,
\end{itemize}
such that they are compatible with the chosen strip-like ends and Floer data, in the sense that
$$\epsilon^*_\zeta K=H_\zeta(t)dt,\quad J(\epsilon_\zeta(s,t))=J_\zeta(t)$$
for each $\zeta\in\Sigma^{\pm}$ and $(s,t)\in Z^{\pm}$.
\end{definition}

For convenience, we call a Floer datum together with a perturbation datum the
{\it analytic data}, and denote it by $\mathbf D_{\mathrm{Fuk}}$.

\subsubsection{Grading of Lagrangian submanifolds}
Let $(\mathbb R^{2n},\omega)$ be the standard symplectic vector space. Denote
by $\Lag_n=\Lag(\mathbb R^{2n},\omega)$ the set of all linear
Lagrangian subspaces. It is known that $\pi_1(\Lag_n)\cong\mathbb Z$ ({\it c.f.} \cite[Theorem 2.31]{McDSal}).
Denote by $\widetilde{\Lag}_n$ the universal covering of $\Lag_n$.

Now suppose $(M^{2n},\omega)$ is a symplectic manifold, then to each $p\in M$ is associated $\Lag(T_pM)$, and one
obtains a fiber bundle, denoted by $\Lag(M)\to M$.

\begin{lemma}\label{fuklemma1}
There exists a covering $\widetilde{\Lag}(M)$ of $\Lag(M)$ such that its restriction to each
fiber is identified with $\widetilde{\Lag}_n\to\Lag_n$ if and only if $c_1(M)=0$.
\end{lemma}

\begin{proof}
See Fukaya \cite[Lemma 2.6]{Fukaya}.
\end{proof}

From now on we fix a covering $\widetilde{\Lag}(M)$ as in the above lemma.
Now suppose $L$ is a Lagrangian submanifold; then there is a canonical section $s$ of the restriction
of $\Lag(M)$ to $L$, which is given by $s(p)=T_pL\subset T_pM$.

\begin{definition}[Grading of Lagrangian Submanifolds]
A {\it graded Lagrangian submanifold} of $(M,\widetilde{\Lag}(M))$ is an oriented Lagrangian submanifold
$L$ and a lift of $s$ to $\tilde s: L\to\widetilde{\Lag}(M)$. The lifting $\tilde s$ is
called the {\it grading} of $L$, and denote $L$ with $\tilde s$ by $\tilde L$.
\end{definition}

The grading of a Lagrangian submanifold is related to its Maslov class as follows:
Let
$\phi:(D^2,\partial D^2)\to (M,L)$ be a map representing
$\pi_2(M,L)$. For each $t\in\partial D^2$ we have a Lagrangian
subspace $T_{\phi(t)}L\subset T_{\phi(t)}M$, which gives a map
$S^1\to\Lag_n$. It determines an element in
$\pi_1(\Lag_n)\cong\mathbb Z$, which is called the {\it Maslov
index} of $[\phi]$, and is denoted by $\mu([\phi])$.
Under the assumption that $c_1(M)=0$, $\mu$ can be
extended to $\pi_1(L)$ as follows: By Lemma \ref{fuklemma1} there exist
a lift $\widetilde{\Lag}(M)\to M$. Let $\gamma:S^1\to L$ be a representation
of an element of $\pi_1(L)$. Define a map
$\gamma^+:S^1\to \Lag(M)$ by
$$\gamma^+(t)=T_{\gamma(t)}L\in\Lag(T_pM).$$
Since $\widetilde{\Lag}(M)\to\Lag(M)$ is a covering, we have a lift
$\tilde\gamma^+:[0,1]\to\widetilde{\Lag}(M)$ of $\gamma^+$.
By the fact that $\widetilde{\Lag}_n/\Z=\Lag_n$ there exists $\ov\mu(\gamma)\in\Z$ such that
$$\ov\mu(\gamma)\cdot\tilde\gamma^+(0)=\tilde\gamma^+(1).$$
The map
$\ov\mu:\pi_1(L)\to\Z$ is called the {\it Maslov class} of $L$. We have:

\begin{lemma}
Suppose $c_1(M)=0$, then there exists a lift $\tilde s$ of $s:L\to\Lag(M)$ if and only if the Maslov class
$\ov\mu:\pi_1(L)\to\Z$ is zero.
\end{lemma}

\begin{proof}
See Fukaya \cite[Lemma 2.14]{Fukaya}.
\end{proof}

\subsubsection{Definition of a Floer cochain}

Suppose $\tilde L_1,\tilde L_2$ intersect transversally and $p\in\tilde L_1\cap \tilde L_2$;
we next define a grading $\eta_{\tilde L_1,\tilde L_2}(p)$ for $p$.
Let $$Y:=D^2\cup\{x+y\sqrt{-1} | x\ge 0, y\in[-1,1]\}\subset \mathbb C.$$
The boundary $\partial Y$ is identified with $\mathbb R$ where $\infty-\sqrt{-1}$ corresponds to $-\infty$
and $\infty+\sqrt{-1}$ corresponds to $+\infty$. Define a path $\tilde l:\mathbb R\to\widetilde{\Lag}(T_pM)$
such that $\tilde l(-\infty)=\tilde s_1(p)$ and $\tilde l(\infty)=\tilde s_2(p)$, where
$\tilde s_1$ and $\tilde s_2$ are the gradings of the Lagrangian submanifolds $L_1$ and $L_2$.
Assume that $\tilde l(t)$ is locally
constant if $|t|>T$.

\begin{lemma}Let $l=\pi\circ\tilde l$ and $W^{1,k}(Y, T_pM; l)
:=\{u\in W^{1,k}(Y,T_pM)|u(x)\in l(x)~\mbox{if }~x\in\partial Y\cong\mathbb R\}$.
Then
\begin{equation}\label{dbar_op}
\ov\partial:W^{1,k}(Y,T_pM;l)\to W^{0,k}(Y,T_pM\otimes\Lambda^{0,1})
\end{equation}
is a Fredholm operator.
\end{lemma}

\begin{proof}See Fukaya \cite[Lemma 3.9]{Fukaya}.\end{proof}

\begin{definition}[Floer Cochain]
If $\tilde L_1,\tilde L_2$ intersect transversally, then
the {\it grading} $\eta_{\tilde L_1,\tilde L_2}(p)$ of $p$ is defined to be the index of $\ov\partial$ in above lemma.
More generally, for two arbitrary $\tilde L_1,\tilde L_2$ with analytic data, let
\begin{equation}\label{defHom}\Hom(\tilde L_1,\tilde L_2):=\mbox{Span}\{y:[0,1]\to M|y(0)\in L_1, y(1)\in L_2,\mbox{ and }dy/dt=X(t, y(t))\}.
\end{equation}
where the grading of $y$, when viewed as the intersection point $p$ of $\tilde L_1,\phi(\tilde L_2)$,
is defined to be the grading $\eta_{\tilde L_1,\phi(\tilde L_2)}(p)$.
An element in $\Hom(\tilde L_1,\tilde L_2)$ is called
a {\it Floer cochain} of $\tilde L_1$ and $\tilde L_2$, and $y$ is sometimes
called a {\it Hamiltonian chord}.
\end{definition}


\begin{lemma}\label{qqstar}
If $\tilde L_1,\tilde L_2$ intersect transversally, then
one may choose $H_{L_0,L_1}$ to be zero,
and $p\in\Hom(\tilde L_1,\tilde L_2)$ implies
the same $p$ lies in $\Hom(\tilde L_2,\tilde L_1)$;
to distinguish, we write $p^*\in\Hom(\tilde L_2,\tilde L_1)$. We have
$$\eta_{\tilde L_0,\tilde L_1}(p)+\eta_{\tilde L_1,\tilde L_0}(p^*)={\dim M}/{2}.$$
\end{lemma}

\begin{proof}
See Fukaya \cite[Lemma 2.27]{Fukaya}.
\end{proof}

\subsubsection{Moduli space of pseudo-holomorphic disks}

Take a pointed-boundary disk $S$ with Lagrangian labels. Equip it with strip-like ends,
Floer data $(H_\zeta,J_\zeta)$ for each point at infinity, and a compatible
perturbation datem $(K,J)$. $K$ determines a vector-field-valued 1-form $Y\in\Omega^1(S,C^\infty(T))$:
for each $\xi\in TS$, $Y(\xi)$ is the Hamiltonian vector field of $K(\xi)$. The
{\it inhomogeneous pseudo-holomorphic map} equation for $u\in C^\infty(S,M)$
is
\begin{equation}\label{IHE}
\left\{
\begin{array}{l}
Du(z)+J(z,u)\circ Du(z)\circ I_S=Y(z,u)+J(z,u)\circ Y(z,u)\circ I_S,\\
u(C)\subset L_C\quad\mbox{for all}\quad C\subset\partial S,
\end{array}
\right.
\end{equation}
where $I_S$ is the complex structure on $S$.

By varying the complex structures on $S$ (we require that at infinity the complex structures is fixed), one obtains
a universal family of pointed-boundary disks with strip-like ends, equipped with Lagrangian labels. For such a family,
one may choose a family of consistent perturbation data (for the existence see Seidel~\cite[\S9i]{Seidel}).
Now suppose
$a_1\in\Hom(\tilde L_1,\tilde L_2),\cdots,a_n\in\Hom(\tilde L_n,\tilde L_{n+1}),a_{n+1}\in\Hom(\tilde L_1,\tilde L_{n+1})$.
Let
$${\mathcal M}(a_1,a_2,\cdots,a_{n+1}):=\left\{u\in C^{\infty}(S,M)\left|
\begin{array}{l}
u\,\mbox{satisfies~(\ref{IHE})\,and the strip-like ends}\\
\mbox{converge to}\, a_1,\cdots,a_{n+1},\,\mbox{respectively}
\end{array}\right.\right\}$$
be the moduli space of solutions to (\ref{IHE}).

\subsubsection{Compactification and orientation of the moduli spaces}

\begin{theorem}
${\mathcal M}$ admits a natural compactification and orientation.
\end{theorem}

\begin{proof}
See Seidel~\cite[\S9l]{Seidel}.\end{proof}

The compactification of $\mathcal M(a_1,a_2,\cdots, a_{n+1})$ is a smooth stratified space (manifold with corners),
where the corners consists of all possible pseudo-holomorphic disks with ``bubbling-off'' disks.
Its codimension one strata consists of
\begin{eqnarray}\label{Fukaya-boundary}
\bigcup_{1\le i<j\le n+1}\bigcup_{b\in \Hom(\tilde{L}_i, \tilde{L}_{j-1})}
\mathcal{M}(b,a_{i},\cdots,a_{j-1})\times\mathcal{M}
(a_{1},\cdots, a_{i-1}, b, a_{j},\cdots, a_{n+1}).
\end{eqnarray}
The orientation is signed the following way: for each $a_i$,
let $o_{a_i}$ be the determinant bundle $\det\ov\partial$ of equation (\ref{dbar_op});
then the orientation bundle of $\mathcal M(a_{n+1},a_1,\cdots, a_n)$ is
$$o_{a_{n+1}}\otimes o_{a_1}^{-}\otimes\cdots\otimes o_{a_n}^-,$$
where $o_{a_i}^-$ is the dual bundle of $o_{a_i}$.
Note that, $\mathcal M(a_{n+1},a_1,\cdots, a_n)$
and $\mathcal M(a_1,\cdots, a_n,a_{n+1})$ count the same set of pseudo-holomorphic disks,
however,
their orientations agree if and only if
$|a_{n+1}|(|a_1|+\cdots +|a_n|)$ is even.

\subsubsection{Construction of the Fukaya category}

\begin{theorem}
Suppose $M$ is an exact symplectic manifold with $c_1(M)=0$,
and possibly with contact type boundary.
Suppose $\tilde L_1,\tilde L_2,\cdots,\tilde L_{n+1}$ are admissible graded
Lagrangian submanifolds, and $a_i\in\Hom(\tilde L_i,\tilde L_{i+1})$,
$i=1,2,\cdots, n$. Define
$$\begin{array}{cccl}
\ov m_n:&\ov{\Hom}(\tilde L_1,\tilde L_2)\otimes\cdots\otimes\ov{\Hom}
(\tilde L_n,\tilde L_{n+1})&\longrightarrow&
\ov{\Hom}(\tilde L_1,\tilde L_{n+1})\\
&(\ov a_1,\ov a_2,\cdots,\ov a_n)&\longmapsto&\displaystyle\sum_{a\in\Hom(\tilde L_1,\tilde L_{n+1})}\#
{\mathcal M}(a,a_1, \cdots, a_n)\cdot \ov a,
\end{array}
$$
for $n=1,2,\cdots$.
Then the set of admissible Lagrangian submanifolds
and the Floer cochain complex among them together with
$\{\ov m_n\}$ defined above form an A$_\infty$ category,
called the {\it Fukaya category} of $M$, and is denoted by
$\mathcal Fuk(M)$.
\end{theorem}

This is proved in \cite[Chapter II]{Seidel} for the exact case and in \cite{Fukaya} for the
general case; we will not repeat it. The (ir)relevance of the construction to the choice
of the analytic data is also completely discussed in \cite[\S12]{Seidel}. Such a
technical problem will also appear in our case when counting the pseudo-holomorphic
disks with punctures. However, all Seidel's argument can be applied to our case,
and we will not address this issue in current paper.

\subsubsection{Cyclicity and a strengthening of the analytic data}\label{subsect_cycl}

From Seidel's original definition, one sees that:
\begin{itemize}
\item if $L_0,L_1$ intersect transversally, then one may choose $H_{L_0,L_1}=0$,
and up to a degree shifting $\Hom(\tilde L_0,\tilde L_1)\cong\Hom(\tilde L_1,\tilde L_0)$
(see Lemma \ref{qqstar}); however,
\item if $L_0, L_1$ do not
intersect transversally, then $H_{L_0, L_1}$ does not vanish, and therefore
$\Hom(\tilde L_0,\tilde L_1)$ is by no means
the same as $\Hom(\tilde L_1,\tilde L_0)$.
\end{itemize}
To overcome this inconsistency, {\it i.e.}
to make Lemma \ref{qqstar} hold even for non transversal pair of Lagrangian submanifolds, we
make the following additional condition for the Floer data:

\begin{definition}[Modified Floer Datum]
For each {\it ordered} pair of Lagrangian submanifolds
$L_0, L_1\subset M$, a {\it Floer datum} consists of
$H_{L_0, L_1}\in C^\infty([0,1], \mathscr H)$ and $J\in C^\infty([0,1],\mathscr J)$, besides the requirement of Definition \ref{F-data}, satisfying the
following additional properties:
\begin{enumerate}
\item for the opposite ordered pair $(L_1, L_0)$, $H_{L_1, L_0}(t)=-H_{L_0, L_1}(1-t)$;
\item if $X$ is the time-dependent Hamiltonian vector field of $H_{L_0, L_1}$ and $\phi_X$
its flow, then $\phi^{1/2}_{X(t)}(L_0)$ intersects $\phi^{1/2}_{-X(1-t)}(L_1)$ transversally. (Note $-X(1-t)$
is the Hamiltonian vector field of $H_{L_1, L_0}$.)
\end{enumerate}
\end{definition}

The perturbation data will be changed accordingly. With such a modification, one sees that
in the case when $L_0$ and $L_1$ do not intersect transversally,
the generators of $\Hom(\tilde L_0,\tilde L_1)$ and $\Hom(\tilde L_1,\tilde L_0)$
may both be identified with the intersection points of $\phi^{1/2}_X(L_0)$ and $\phi^{1/2}_{-X}(L_1)$,
and the degrees at each point add up to $n$.
An application of this is the Lie bialgebra structure on the cyclic complex of
the Fukaya category, which we discuss in the next subsection.

\subsection{The Lie bialgebra structure}\label{sect_liebialg}

From now on, we graded the Floer cochains negatively. Such a convention
is usually adopted in algebraic topology when studying
Hochschild/cyclic homology of the cochain complex of topological spaces.

\begin{definition}[Lie Bialgebra]\label{def_Liebialg} Let $L$ be a (possibly graded)
$\mathbb K$-space. A Lie bialgebra on $L$ is the triple $(L, [\,,\,],\delta)$ such that
\begin{itemize}\item $(L, [\,,\,])$ is a Lie algebra;
\item $(L, \delta)$ is a Lie coalgebra;
\item The Lie algebra and coalgebra satisfy the following identity, called the {\it Drinfeld compatibility}:
$$\delta[a,b]=\sum_{(a)}((-1)^{|a''||b|}[a',b]\otimes a''+a'\otimes [a'', b])+\sum_{(b)}([a,b']\otimes b''+(-1)^{|a||b'|}b'\otimes[a,b'']),$$
for all $a,b\in L$, where we write $\delta(a)=\sum_{(a)}a'\otimes a''$ and $\delta(b)=\sum_{(b)}b'\otimes b''$.
\end{itemize}
If moreover, $[\,,\,]\circ\delta(a)\equiv 0$, for all $a\in L$,
$(L,[\,,\,],\delta)$ is said to be {\it involutive}. If the Lie
bracket has degree $k$ and the Lie cobracket has degree $l$, denote
the Lie bialgebra with degree $(l,k)$.
\end{definition}

\begin{theorem}[Lie Bialgebra of The Fukaya Category]\label{Liebi_Fuk}
Let $M^{2n}$ be an exact symplectic manifold (possibly with contact type boundary) with $c_1(M)=0$.
Grade the Floer cochain complex negatively.
Then the cyclic cochain complex of the Fukaya category $\mathcal Fuk(M)$ of $M$ has the structure
of a differential involutive Lie bialgebra of degree $(2-n, 2-n)$.
\end{theorem}

The proof of this theorem consists of the rest of the subsection. Before going to
the details, we would like to say several words about the degrees.
We say a graded vector space $V$ is a Lie algebra of degree $n$ if $V[-n]$ is a
graded Lie algebra in the usual sense.
Similar convention applies to Lie bialgebras with a bi-degree $(m,n)$. A technical issue here
is the correctness of Drinfeld compatibility and involutivity.
In our case of the Lie bialgebra of degree $(2-n,2-n)$, if we shift the vector space down
by $2-n$, then the Lie bracket has degree zero and the Lie cobracket has degree $4-2n$, which is even.
Therefore, equations for the Drinfeld compatibility and involutivity in this case
is the same as in the usual case.

\begin{lemma}[Lie Algebra]
Denote by $\mathcal Fuk(M)$ the Fukaya category of $M$.
Define $$[\,,\,]:\Cycl^*(\mathcal Fuk(M))\otimes\Cycl^*(\mathcal Fuk(M))\to\Cycl^*(\mathcal Fuk(M))$$
by
\begin{eqnarray*}[f,g](\overline a_1,\overline a_2,\cdots,\overline a_{n})
&:=&\sum_{i< j}\sum_{p\in\Hom(\tilde L_j,\tilde L_i)}\pm
f(\overline a_i,\cdots,\overline a_{j-1},\overline p)\cdot g(\overline p^*,
\overline a_{j},\cdots, \overline a_{n},\overline a_1, \cdots, \overline a_{i-1})\\
&-&\sum_{i< j}\sum_{p\in\Hom(\tilde L_j, \tilde L_i)}\pm g(\overline a_i,\cdots,
\overline a_{j-1},\overline p)\cdot f(\overline p^*, \overline a_{j},\cdots, \overline a_{n}, \overline a_1,\cdots, \overline a_{i-1}).
\end{eqnarray*}
Then $(\Cycl^*(\mathcal Fuk(M)),[\,,\,],b)$ forms
a differential graded Lie algebra of degree $2-n$.
\end{lemma}

Pictorially, the bracket is defined as in the following picture (Figure \ref{fig1}):
\begin{figure}[h]
$$
\setlength{\unitlength}{.9cm}
\begin{picture}(15,4)
\thicklines
\put(2.5,1.5){\circle{3.6}}
\put(2.2, 1.5){$[f,g]$}

\put(4.5, 1.5){$\tilde L_1$}
\put(4.12,2.3){\circle*{.1}}
\put(4.3,2.3){$a_1$}

\put(3, 3.4){$\tilde L_2$}
\put(1.9,3.2){\circle*{.1}}
\put(1.7,3.4){$a_2$}

\put(.92,2.4){\circle*{.1}}
\put(.92,.62){\circle*{.1}}

\put(4.3,.5){$a_n$}
\put(4.15,.8){\circle*{.1}}
\put(3.5,-0.4){$\tilde L_n$}

\put(5.2,1.5){$=\sum\limits_{p}$}

\put(7,1.5){\circle{1.6}}
\put(8.6,1.5){\circle{1.6}}
\put(7.8,1.5){\circle*{.1}}
\put(7.5,1.4){$p$}
\put(7.9,1.4){$p^*$}
\put(8.6,1.4){$g$}
\put(6.8,1.4){$f$}

\put(7.2,2.4){$a_i$}\put(7.25,2.25){\circle*{.1}}
\put(7.8,2.4){$a_{i-1}$}\put(8.3,2.25){\circle*{.1}}
\put(6.9,.4){$a_{j-1}$}\put(7.25,.75){\circle*{.1}}
\put(8,.4){$a_{j}$}\put(8.2,.8){\circle*{.1}}
\put(9.2,2.05){\circle*{.1}}
\put(9.2,.95){\circle*{.1}}
\put(6.2,1.5){\circle*{.1}}

\put(9.6,1.5){$-\sum\limits_p$}

\put(11.4,1.5){\circle{1.6}}
\put(10.6,1.5){\circle*{.1}}
\put(13,1.5){\circle{1.6}}
\put(12.2,1.5){\circle*{.1}}
\put(11.9,1.4){$p$}
\put(12.3,1.4){$p^*$}
\put(11.2,1.4){$g$}
\put(13,1.4){$f$}

\put(11.6,2.4){$a_i$}\put(11.65,2.25){\circle*{.1}}
\put(12.2,2.4){$a_{i-1}$}\put(12.7,2.25){\circle*{.1}}
\put(11.3,.4){$a_{j-1}$}\put(11.65,.75){\circle*{.1}}
\put(12.4,.4){$a_{j}$}\put(12.6,.8){\circle*{.1}}
\put(13.6,2.05){\circle*{.1}}
\put(13.6,.95){\circle*{.1}}
\end{picture}
$$
\caption{Definition of the Lie bracket for two cyclic cochains}\label{fig1}
\end{figure}
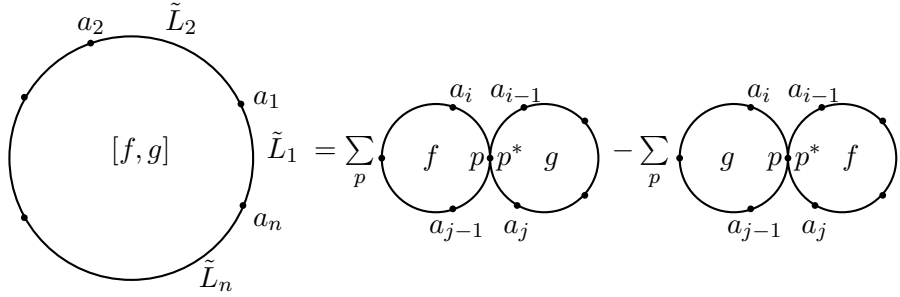
in the picture, the left side of the equality is
the value of $[f,g]$ on $(a_1,\cdots,a_n)$, and the right side
of the equality is summarized over all possibilities of the product of the value of $f$ on $(a_i,\cdots,a_{j-1},p)$
with the value of $g$ on $(p^*,a_j,\cdots,a_{i-1})$.

\begin{proof}
First, we show $[\,,\,]$ has degree $2-n$. Note that $(\overline a_1,\cdots,\overline a_n)$ has degree
$$|a_1|+\cdots+|a_n|+n,$$
(recall that we grade $a_i\in \Hom(\tilde L_i, \tilde L_{i+1})$ negatively). And the sum of the degrees of
$(\overline a_i,\cdots,\overline a_{j-1},\overline p)$ and $(\overline p^*, \overline a_{j},\cdots, \overline a_{n}, \cdots, \overline a_{i-1})$
is
\begin{eqnarray*}
&&|a_1|+\cdots+|a_n|+|p|+|p^*|+(n+2)\\
&=&|a_1|+\cdots+|a_n|+|p|+(-n-|p|)+(n+2)\quad(\mbox{recall that}~|p^*|+|p|=-n)\\
&=&|a_1|+\cdots+|a_n|+n+(2-n).
\end{eqnarray*}
The difference of these two degrees is exactly $n-2$.
Going to the cyclic cochain ({\it i.e.} the dual space) level, the degree of the bracket
$[\,,\,]$ becomes $2-n$.

Second, we show $[\,,\,]$ is graded skew-symmetric. Observe that in the definition of $[\,,\,]$, if we switch $f$ and $g$,
we get exactly the opposite sign (ignoring the intrinsic signs that come from the Koszul convention).

Third, we show the Jacobi identity: With Figure~\ref{fig1} in mind,
the value of $[[f,g],h]$ on $(\ov a_1,\cdots,\ov a_n)$ has four terms which can be pictorially
represented by the following picture
$$
\setlength{\unitlength}{.8cm}
\begin{picture}(11,2)
\put(.5,.5){\circle{1}}
\put(.5,1.5){\circle{1}}
\put(1.5,.5){\circle{1}}
\put(.4,.4){$g$}
\put(.4,1.4){$f$}
\put(1.4,.4){$h$}

\put(2.3,1){$-$}

\put(3.5,.5){\circle{1}}
\put(3.5,1.5){\circle{1}}
\put(4.5,.5){\circle{1}}
\put(3.4,.4){$f$}
\put(3.4,1.4){$g$}
\put(4.4,.4){$h$}

\put(5.3,1){$-$}

\put(6.5,1.5){\circle{1}}
\put(7.5,1.5){\circle{1}}
\put(7.5,.5){\circle{1}}
\put(6.4,1.4){$h$}
\put(7.4,1.4){$f$}
\put(7.4,.4){$g$}

\put(8.3,1){$+$}

\put(9.5,1.5){\circle{1}}
\put(10.5,1.5){\circle{1}}
\put(10.5,.5){\circle{1}}
\put(9.4,1.4){$h$}
\put(10.4,1.4){$g$}
\put(10.4,.4){$f$}
\end{picture}
$$
Similarly, $[[g,h],f]$ and $[[h,f],g]$
are represented by the following picture:
$$
\setlength{\unitlength}{.8cm}
\begin{picture}(11,4.5)
\put(.5,.5){\circle{1}}
\put(.5,1.5){\circle{1}}
\put(1.5,.5){\circle{1}}
\put(.4,.4){$h$}
\put(.4,1.4){$g$}
\put(1.4,.4){$f$}

\put(2.3,1){$-$}

\put(3.5,.5){\circle{1}}
\put(3.5,1.5){\circle{1}}
\put(4.5,.5){\circle{1}}
\put(3.4,.4){$g$}
\put(3.4,1.4){$h$}
\put(4.4,.4){$f$}

\put(5.3,1){$-$}

\put(6.5,1.5){\circle{1}}
\put(7.5,1.5){\circle{1}}
\put(7.5,.5){\circle{1}}
\put(6.4,1.4){$f$}
\put(7.4,1.4){$g$}
\put(7.4,.4){$h$}

\put(8.3,1){$+$}

\put(9.5,1.5){\circle{1}}
\put(10.5,1.5){\circle{1}}
\put(10.5,.5){\circle{1}}
\put(9.4,1.4){$f$}
\put(10.4,1.4){$h$}
\put(10.4,.4){$g$}

\put(.5,3){\circle{1}}
\put(.5,4){\circle{1}}
\put(1.5,3){\circle{1}}
\put(.4,2.9){$f$}
\put(.4,3.9){$h$}
\put(1.4,2.9){$g$}

\put(2.3,3.5){$-$}

\put(3.5,3){\circle{1}}
\put(3.5,4){\circle{1}}
\put(4.5,3){\circle{1}}
\put(3.4,2.9){$h$}
\put(3.4,3.9){$f$}
\put(4.4,2.9){$g$}

\put(5.3,3.5){$-$}

\put(6.5,4){\circle{1}}
\put(7.5,4){\circle{1}}
\put(7.5,3){\circle{1}}
\put(6.4,3.9){$g$}
\put(7.4,3.9){$h$}
\put(7.4,2.9){$f$}

\put(8.3,3.5){$+$}

\put(9.5,4){\circle{1}}
\put(10.5,4){\circle{1}}
\put(10.5,3){\circle{1}}
\put(9.4,3.9){$g$}
\put(10.4,3.9){$f$}
\put(10.4,2.9){$h$}\end{picture}
$$
The sum is identically zero, which proves the (graded) Jacobi identity.

Finally, we show that the bracket commutes with the boundary:
\begin{eqnarray}&&(b[f,g])(\ov a_1,\ov a_2,\cdots,\ov a_{n})\nonumber\\
&=&[f,g](b(\ov a_1,\ov a_2,\cdots, \ov a_{n}))\nonumber\\
&=&[f,g]\big(\sum_{i}\sum_{k}\pm(\ov a_1,\cdots,\overline m_k(\ov a_{i},\cdots,\ov a_{i+k-1}),\cdots,\ov a_n)\big)\label{Eq_1'}\\
&+&[f,g]\big(\sum_{j}\sum_{k}\pm(\overline m_k(\ov
a_{n-j},\cdots,\ov a_n,\ov a_1,\cdots,\ov a_{i}),\cdots,\ov
a_{n-j-1})\big),\label{Eq_1''}
\end{eqnarray}
while
\begin{eqnarray}
&&([bf,g]+(-1)^{|f|}[f, bg])(\ov a_1,\ov a_2,\cdots,\ov a_{n})\nonumber\\
&=&\sum_{i< j}\sum_{p}\pm f(b(\ov a_i,\cdots,\ov a_{j-1},\ov p))\cdot g(\ov p^*,\ov a_{j},\cdots,\ov a_{i-1})\label{Eq_2}\\
&-&\sum_{i< j}\sum_{p}\pm g(\ov a_i,\cdots,\ov a_{j-1},\ov p)\cdot f(b(\ov p^*,\ov a_{j},\cdots,\ov a_{i-1}))\label{Eq_22}\\
&+&\sum_{i< j}\sum_{p}\pm f(\ov a_i,\cdots,\ov a_{j-1}, \ov p)\cdot g(b(\ov p^*,\ov a_{j},\cdots,\ov a_{i-1}))\label{Eq_3}\\
&-&\sum_{p}\pm g(b(\ov a_i,\cdots,\ov a_{j-1}, \ov p))\cdot f(\ov p^*,\ov a_{j},\cdots,\ov a_{i-1}).\label{Eq_32}
\end{eqnarray}
From the definition of $[\,,\,]$, one sees $(\ref{Eq_2})+(\ref{Eq_22})+(\ref{Eq_3})+(\ref{Eq_32})$
contains more terms than $(\ref{Eq_1'})+(\ref{Eq_1''})$, namely, those
terms involving $\overline m_k$ acting on $\ov p$ and $\ov p^*$. For example, the extra terms coming from
$(\ref{Eq_2})$ are
\begin{equation}
\sum_{p\in\Hom(\tilde L_j,\tilde L_i)}\sum_k\sum_r f(\ov m_r(\ov a_k,\cdots, \ov a_{j-1}, \ov p, \ov a_{i},\cdots,\ov a_l),\ov a_{l+1},\cdots,
\ov a_{k-1})\cdot g(\ov p^*,\ov a_{j},\cdots, \ov a_{i-1})\label{Eq_4}
\end{equation}
and the ones from $(\ref{Eq_3})$ are
\begin{equation}
\sum_{p\in\Hom(\tilde L_j,\tilde L_i)}\sum_k\sum_r f(\ov a_i,\cdots, \ov a_j,\ov p)\cdot
g(\overline m_r(\ov a_k,\cdots, \ov a_{i-1}, \ov p^*,\ov a_{j},\cdots,\ov a_l),\ov a_{l+1},\cdots,\ov a_{k-1}).\label{Eq_5}
\end{equation}
However, these two groups of terms cancel with each other because
\begin{eqnarray*}&&\sum_{p\in\Hom(\tilde L_j,\tilde L_i)}
\overline m_r(\ov a_k, \cdots,\ov a_j,\ov p, \ov a_i, \cdots,\ov a_l)\otimes\ov p^*\\
&=&\sum_{p\in\Hom(\tilde L_j,\tilde L_i)}\sum_{q\in\Hom(\tilde L_{k},\tilde L_{l+1})}
\#\mathcal M(q,a_k, \cdots, a_j,  p, a_i, \cdots, a_l)\ov q\otimes\ov p^*\\
&=&\sum_{q\in\Hom(\tilde L_{k},\tilde L_{l+1})}
\sum_{p\in\Hom(\tilde L_j,\tilde L_i)}\ov q\otimes\# \mathcal M(q,a_k, \cdots, a_j, p, a_i, \cdots, a_l )\ov p^*\\
&\stackrel{\mbox{\scriptsize cyclicity}}{=}
&\sum_{q\in\Hom(\tilde L_{k},\tilde L_{l+1})}\sum_{p\in\Hom(\tilde L_j,\tilde L_i)}
\ov q\otimes\#\mathcal M(p,a_i,\cdots,a_l,q,a_k, \cdots, a_j )\ov p^*\\
&\stackrel{\S\ref{subsect_cycl}}{=}&\sum_{q\in\Hom(\tilde L_{k},\tilde L_{l+1})}\sum_{p^*\in\Hom(\tilde L_i,\tilde L_j)}
\ov q\otimes\#\mathcal M(p^*,a_i,\cdots,a_l,q^*,a_k, \cdots, a_j )\ov p^*\\
&=&\sum_{q\in\Hom(\tilde L_{k},\tilde L_{l+1})}\ov q\otimes \overline m_r(\ov a_i,\cdots,\ov a_l,\ov q^*,\ov a_k,\cdots,\ov a_j).
\end{eqnarray*}
By substituting the above identity into (\ref{Eq_4}) we get exactly (\ref{Eq_5}).
Similarly, the extra terms in (\ref{Eq_22}) and in (\ref{Eq_32}) cancel with each other.
Pictorially, the value of $b[f,g]$ on $(a_1,\cdots, a_n)$ equals
$$
\setlength{\unitlength}{.9cm}
\begin{picture}(12,2.5)
\put(.5,.5){\circle{1}}
\put(.5,1.5){\circle{1}}
\put(1.5,.5){\circle{1}}
\put(.4,.4){$f$}
\put(.1,1.4){$\#\mathcal M$}
\put(1.4,.4){$g$}

\put(2.3,1){$-$}

\put(3.5,.5){\circle{1}}
\put(3.5,1.5){\circle{1}}
\put(4.5,.5){\circle{1}}
\put(3.4,.4){$g$}
\put(3.1,1.4){$\#\mathcal M$}
\put(4.4,.4){$f$}

\put(5.3,1){$+$}

\put(6.5,.5){\circle{1}}
\put(7.5,1.5){\circle{1}}
\put(7.5,.5){\circle{1}}
\put(6.4,.4){$f$}
\put(7.1,1.4){$\#\mathcal M$}
\put(7.4,.4){$g$}

\put(8.3,1){$-$}

\put(9.5,.5){\circle{1}}
\put(10.5,1.5){\circle{1}}
\put(10.5,.5){\circle{1}}
\put(9.4,.4){$g$}
\put(10.1,1.4){$\#\mathcal M$}
\put(10.4,.4){$f$}
\end{picture}
$$
and the value of $[bf,g]+[f,bg]$ on $(a_1,\cdots,a_n)$ not only contains the above four terms, but also
$$
\setlength{\unitlength}{.9cm}
\begin{picture}(9,2.5)
\put(-.3,.5){$-$}
\put(1,.5){\circle{1}}\put(.9,.4){$f$}
\put(2,.5){\circle{1}}\put(1.6,.4){$\#\mathcal M$}
\put(3,.5){\circle{1}}\put(2.9,.4){$g$}
\put(3.8,.5){$+$}
\put(5,.5){\circle{1}}\put(4.9,.4){$g$}
\put(6,.5){\circle{1}}\put(5.6,.4){$\#\mathcal M$}
\put(7,.5){\circle{1}}\put(6.9,.4){$f$}

\put(1,2){\circle{1}}\put(.9,1.9){$f$}
\put(2,2){\circle{1}}\put(1.6,1.9){$\#\mathcal M$}
\put(3,2){\circle{1}}\put(2.9,1.9){$g$}
\put(3.8,2){$-$}
\put(5,2){\circle{1}}\put(4.9,1.9){$g$}
\put(6,2){\circle{1}}\put(5.6,1.9){$\#\mathcal M$}
\put(7,2){\circle{1}}\put(6.9,1.9){$f$}

\end{picture}
$$
which cancel each other within themselves.
\end{proof}

\begin{lemma}[Lie Coalgebra]
Denote by $\mathcal Fuk(M)$ the Fukaya category of $M$.
Define $\Cycl^*(\mathcal Fuk(M))\to
\Cycl^*(\mathcal Fuk(M))\otimes\Cycl^*(\mathcal Fuk(M))$
by
\begin{eqnarray*}
&&(\delta f)(\ov a_1, \ov a_2,\cdots,\ov a_{n})\otimes(\ov b_1,\ov b_2,\cdots, \ov b_{m})\\
&:=&\sum_{i=1}^n\sum_{j=1}^m\sum_{p\in\Hom(\tilde L_i,\tilde L_j)}\pm
f(\ov a_1,\cdots,\ov a_{i-1},\ov p,\ov b_j,\cdots, \ov b_{m},\ov b_1,\cdots,\ov b_{j-1},
\ov p^*,\ov a_i,\cdots, \ov a_{n}).
\end{eqnarray*}
Then $(\Cycl^*(\mathcal Fuk(M)),\delta,b)$ forms a Lie coalgebra of degree $2-n$.
\end{lemma}

Pictorially, the cobracket is defined as follows:
$$
\setlength{\unitlength}{.8cm}
\begin{picture}(8,1.2)

\put(.3,.4){$(\delta f)\Big($}

\put(2,.5){\circle{1}}\put(1.9,.4){$x$}

\put(2.6,.4){$\otimes$}

\put(3.5,.5){\circle{1}}\put(3.4,.4){$y$}

\put(4.1,.4){$\Big)=$}

\put(5,.4){$f\Big($}

\put(6,.5){\arc[20,340]{.5}}\put(5.9,.4){$x$}
\put(6.95,.5){\arc[-160,160]{.5}}\put(6.85,.4){$y$}
\put(6.475,.7){\circle*{.1}}\put(6.375,1){$p$}
\put(6.475,.3){\circle*{.1}}\put(6.375,-.3){$p^*$}

\put(7.5,.4){$\Big)$}
\end{picture}
$$
In the picture, circled $x$ means $(a_1,\cdots,a_n)$, circled $y$ means $(b_1,\cdots,b_m)$,
and circled $xy$ together in the right side means
$(\ov a_1,\cdots,\ov a_{i-1},\ov p,\ov b_j,\cdots, \ov b_{m},\ov b_1,\cdots,\ov b_{j-1},
\ov p^*,\ov a_i,\cdots, \ov a_{n})$.

\begin{proof} From the definition of $\delta$, the following two statements are obvious:\begin{enumerate}
\item[(1)] $\delta f$ is well defined, namely, the value of $\delta f$ is invariant under the cyclic permutations of
$(\ov a_1,\cdots,\ov a_n)$ and $(\ov b_1,\cdots,\ov b_m)$;
\item[(2)] $\delta f$ is (graded) skew-symmetric, namely,
if we switch $(\ov a_1, \cdots,\ov a_n)$ and $(\ov b_1, \cdots,\ov b_m)$,
the sign of the value of $\delta f$ changes.
\end{enumerate}

The co-Jacobi identity can be proved in a similar way to the proof of Jacobi identity.
Let $\tau:x\otimes y\otimes z\mapsto\pm z\otimes x\otimes y$ be the cyclic permutation of three elements,
then $(\tau^2+\tau+id)\circ(id\otimes\delta)\circ\delta f$ has six terms, grouped into three pairs, pictorially as follows:
$$
\setlength{\unitlength}{.8cm}
\begin{picture}(17,2)

\put(1,.5){\arc[20,70]{.5}}\put(.9,.4){$x$}
\put(1,.5){\arc[110,340]{.5}}
\put(1.95,.5){\arc[-160,160]{.5}}\put(1.85,.4){$y$}
\put(1,1.45){\arc[-70,250]{.5}}\put(.9,1.35){$z$}
\put(1.2,.975){\circle*{.1}}
\put(.8,.975){\circle*{.1}}
\put(1.475,.7){\circle*{.1}}
\put(1.475,.3){\circle*{.1}}
\put(2.8,.875){$-$}
\put(4.95,.5){\arc[110,160]{.5}}\put(3.9,.4){$x$}
\put(4.95,.5){\arc[-160,70]{.5}}
\put(4,.5){\arc[20,340]{.5}}\put(4.85,.4){$y$}
\put(4.95,1.45){\arc[-70,250]{.5}}\put(4.85,1.35){$z$}
\put(5.15,.975){\circle*{.1}}
\put(4.75,.975){\circle*{.1}}
\put(4.475,.7){\circle*{.1}}
\put(4.475,.3){\circle*{.1}}

\put(7,.5){\arc[20,70]{.5}}\put(6.9,.4){$z$}
\put(7,.5){\arc[110,340]{.5}}
\put(7.95,.5){\arc[-160,160]{.5}}\put(7.85,.4){$x$}
\put(7,1.45){\arc[-70,250]{.5}}\put(6.9,1.35){$y$}
\put(7.2,.975){\circle*{.1}}
\put(6.8,.975){\circle*{.1}}
\put(7.475,.7){\circle*{.1}}
\put(7.475,.3){\circle*{.1}}
\put(8.8,.875){$-$}
\put(10.95,.5){\arc[110,160]{.5}}\put(9.9,.4){$z$}
\put(10.95,.5){\arc[-160,70]{.5}}
\put(10,.5){\arc[20,340]{.5}}\put(10.85,.4){$x$}
\put(10.95,1.45){\arc[-70,250]{.5}}\put(10.85,1.35){$y$}
\put(11.15,.975){\circle*{.1}}
\put(10.75,.975){\circle*{.1}}
\put(10.475,.7){\circle*{.1}}
\put(10.475,.38){\circle*{.1}}

\put(13,.5){\arc[20,70]{.5}}\put(12.9,.4){$y$}
\put(13,.5){\arc[110,340]{.5}}
\put(13.95,.5){\arc[-160,160]{.5}}\put(13.85,.4){$z$}
\put(13,1.45){\arc[-70,250]{.5}}\put(12.9,1.35){$x$}
\put(13.2, .975){\circle*{.1}}
\put(12.8, .975){\circle*{.1}}
\put(13.475, .7){\circle*{.1}}
\put(13.475, .3){\circle*{.1}}
\put(14.8, .875){$-$}
\put(16.95, .5){\arc[110,160]{.5}}\put(15.9, .4){$y$}
\put(16.95, .5){\arc[-160,70]{.5}}
\put(16, .5){\arc[20,340]{.5}}\put(16.85, .4){$z$}
\put(16.95,1.45){\arc[-70,250]{.5}}\put(16.85,1.35){$x$}
\put(17.15, .975){\circle*{.1}}
\put(16.75, .975){\circle*{.1}}
\put(16.475, .7){\circle*{.1}}
\put(16.475, .3){\circle*{.1}}

\end{picture}
$$
and they cancel with each other. We obtain the co-Jacobi identity.

Next, we show that $b$ respects the cobracket. This is also similar to the Lie case. By definition,
\begin{eqnarray}
& &((b\otimes id\pm id\otimes b)\delta (f))((\ov a_1, \ov a_2,\cdots,\ov a_{n})
\otimes(\ov b_1,\ov b_2,\cdots, \ov b_{m}))\nonumber\\
&=&\delta f(b(\ov a_1, \ov a_2,\cdots,\ov a_{n})\otimes(\ov b_1,\ov b_2,\cdots, \ov b_{m})
\pm(\ov a_1, \ov a_2,\cdots,\ov a_{n})\otimes b(\ov b_1,\ov b_2,\cdots, \ov b_{m})),\label{codiff1}
\end{eqnarray}
while
\begin{eqnarray}
&&\delta(b(f))((\ov a_1, \ov a_2,\cdots,\ov a_{n})\otimes(\ov b_1,\ov b_2,\cdots, \ov b_{m}))\nonumber\\
&=&\sum_{i=1}^n\sum_{j=1}^m\sum_{p}\pm b(f)
(\ov a_1,\cdots,\ov a_{i-1},\ov p,\ov b_j,\cdots, \ov b_{m},\ov b_1,\cdots,\ov b_{j-1},
\ov p^*,\ov a_i,\cdots, \ov a_{n})\label{codiff3}
\end{eqnarray}
Compared with $(\ref{codiff1})$, $(\ref{codiff3})$ has extra terms
\begin{eqnarray}
&&\sum_{p}f(\ov a_1,\cdots, \overline
m_r(\ov a_{k},\cdots, \ov p,\cdots,\ov b_{l-1}),\ov b_l,\cdots, \ov p^*,\cdots, \ov a_n)\label{cob1}\\
&+&\sum_{p}f(\ov a_1,\cdots,\ov e_p,\cdots,\ov b_{j-1},
\ov m_r(\ov b_j,\cdots, \ov p^*,\cdots,\ov a_l),\cdots,\ov a_n)\label{cob2}\\
&+&\sum_{p}f(\ov a_1,\cdots,\ov m_r(\ov a_k,\cdots,\ov p,\ov b_j,\cdots,
\ov p^*,\cdots, \ov a_{l-1}),\ov a_{l},\cdots,\ov a_n),\label{cob3}
\end{eqnarray}
However, a similar argument as in the Lie case, $(\ref{cob1})+(\ref{cob2})$ vanishes, and the terms in
$(\ref{cob3})$ come in pair (counting $p$ and $p^*$), which
cancel within themselves.
A pictorial proof is also easy, and is left to the interested reader.
This proves that the differential commutes with the cobracket.
\end{proof}

\subsubsection{Proof of the Drinfeld compatibility}
Pictorially, the value of $\delta\circ [f,g]$ on $(a_1,\cdots, a_n)\otimes(b_1,\cdots,b_m)$
is represented by
$$
\setlength{\unitlength}{.8cm}
\begin{picture}(5,1.2)
\put(.3,.4){$[f,g]\Big($}

\put(2,.5){\arc[20,340]{.5}}
\put(2.95,.5){\arc[-160,160]{.5}}
\put(2.475,.7){\circle*{.1}}
\put(2.475,.3){\circle*{.1}}

\put(3.5,.4){$\Big),$}
\end{picture}
$$
which, by definition of $[\,,\,]$, is equal to
$$
\setlength{\unitlength}{.8cm}
\begin{picture}(13,2.3)

\put(0,.9){$-$}

\put(1,.5){\arc[20,340]{.5}}\put(.9,.4){$f$}
\put(1.95,.5){\arc[-160,160]{.5}}
\put(1,1.5){\circle{1}}\put(.9,1.4){$g$}
\put(1.475,.7){\circle*{.1}}
\put(1.475,.3){\circle*{.1}}
\put(2.8,.875){$+$}
\put(4.95,.5){\arc[-160,160]{.5}}
\put(4,.5){\arc[20,340]{.5}}\put(4.85,.4){$f$}
\put(4.95,1.5){\circle{1}}\put(4.85,1.4){$g$}
\put(4.475,.7){\circle*{.1}}
\put(4.475,.3){\circle*{.1}}

\put(6.5,.875){$+$}

\put(8,.5){\arc[20,340]{.5}}\put(7.9,.4){$g$}
\put(8.95,.5){\arc[-160,160]{.5}}
\put(8,1.5){\circle{1}}\put(7.9,1.4){$f$}
\put(8.475,.7){\circle*{.1}}
\put(8.475,.3){\circle*{.1}}
\put(9.8,.875){$-$}
\put(11.95,.5){\arc[-160,160]{.5}}
\put(11,.5){\arc[20,340]{.5}}\put(11.85,.4){$g$}
\put(11.95,1.5){\circle{1}}\put(11.85,1.4){$f$}
\put(11.475,.7){\circle*{.1}}
\put(11.475,.3){\circle*{.1}}

\end{picture}
$$
The left two terms give $[\delta f, g]$ and the right two terms give
$[f,\delta g]$, and we obtain the Drinfeld compatibility.

\subsubsection{Proof of the involutivity}
Suppose $\delta(f)=\sum f'\otimes f''$, then the values of $[f',f'']$ are represent
by
$$
\setlength{\unitlength}{.8cm}
\begin{picture}(11,1.2)

\put(.5,.5){\circle{1}}\put(.4,.4){$f'$}
\put(1.5,.5){\circle{1}}\put(1.4,.4){$f''$}
\put(1,.5){\circle*{.1}}

\put(2.4,.4){$-$}

\put(3.5,.5){\circle{1}}\put(3.4,.4){$f''$}
\put(4.5,.5){\circle{1}}\put(4.4,.4){$f'$}
\put(4,.5){\circle*{.1}}

\put(5.4,.4){$=$}

\put(6.5,.5){\circle{1}}\put(6.4,.4){$f'$}
\put(8,.5){\circle{1}}\put(7.9,.4){$f''$}
\put(7,.5){\circle*{.1}}
\put(7.1,.4){$\times$}
\put(7.5,.5){\circle*{.1}}

\put(8.6,.4){$-$}

\put(9.5,.5){\circle{1}}\put(9.4,.4){$f''$}
\put(11,.5){\circle{1}}\put(10.9,.4){$f'$}
\put(10,.5){\circle*{.1}}
\put(10.1,.4){$\times$}
\put(10.5,.5){\circle*{.1}}\end{picture}
$$
which, by definition of $\delta$, is equal to
$$
\setlength{\unitlength}{.8cm}
\begin{picture}(5,1.5)

\put(.4,.4){$f$}
\put(.5,.5){\arc[60,390]{.5}}
\put(1,.5){\circle*{.1}}
\put(1.2,1.2){\arc[-120,210]{.5}}
\put(.7,1.2){\circle*{.1}}
\put(.78,.92){\circle*{.1}}
\put(.92,.78){\circle*{.1}}

\put(2.3,.8){$-$}

\put(3.4,.4){$f$}
\put(3.5,.5){\arc[60,390]{.5}}
\put(4,.5){\circle*{.1}}
\put(4.2,1.2){\arc[-120,210]{.5}}
\put(3.7,1.2){\circle*{.1}}
\put(3.78,.92){\circle*{.1}}
\put(3.92,.78){\circle*{.1}}
\end{picture}
$$
which is identically zero. For the convenience of readers, 
let us write down the formulas.
By definition, for any $f\in\Cycl^*(\Fuk(M))$, 
\begin{eqnarray*}
&&([\,,\,]\circ \delta(f))(\ov a_1,\ov a_2,\cdots,\ov a_n)\\
&=&\delta(f)\Big(\sum_{i<j}\sum_p\pm(\ov a_i,\cdots,\ov a_{j-1},\ov p)\otimes(\ov a_1,\cdots,
\ov a_{i-1},\ov p^*,\ov a_j,\cdots,\ov a_n)\Big)\\
&-&\delta(f)\Big(\sum_{i<j}\sum_p\pm(\ov a_1,\cdots, \ov a_{i-1},\ov p^*,\ov a_j,\cdots,\ov a_n)\otimes
(\ov a_i,\cdots,\ov a_{j-1},\ov p)\Big).
\end{eqnarray*}
The right hand side of above equality should vanish
because the value of the first half, which is the value of $f$ at 
\begin{eqnarray*}
&&\sum_{i<j}\sum_{i\le k\le j,l<i}\sum_{p,q}
\pm(
\ov a_i,\cdots,\ov a_{k-1},\ov q,\ov a_l,\cdots,\ov a_{i-1},\ov p^*,\ov a_j,
\cdots,\ov a_n,\ov a_1,\cdots,\ov a_{l-1},\ov q^*,\ov a_k,\cdots,\ov a_{j-1},\ov p)\\
&+&\sum_{i<j}\sum_{i\le k\le j,j<l}\sum_{p,q}
\pm(
\ov a_i,\cdots,\ov a_{k-1},\ov q,\ov a_l,\cdots,\ov a_n,\ov a_1,\cdots,\ov a_{i-1},\ov p^*,\ov a_j,
\cdots,\ov a_{l-1},\ov q^*,\ov a_k,\cdots,\ov a_{j-1},\ov p),\end{eqnarray*}
is the same as the value of $f$ at the second half up to a cyclic order. This proves the involutivity.

\section{Linearized contact homology}\label{sect_lin}

The contact homology of a contact manifold was first introduced in symplectic field theory
by Eliashberg-Givental-Hofer (\cite{EGH}) in late 1990s.
Its linearized version, the {\it linearized contact homology}, can be found
in \cite{CL} and \cite{BO}. Let us recall its definition.

\subsection{Several concepts in contact geometry}

Let $W$ be a manifold of dimension $2n-1$. 
A {\it contact form} on $W$
is a 1-form $\lambda$ such that $\lambda\wedge(d\lambda)^{n-1}$ is a volume form on $W$
(here we only consider co-orientable contact manifolds).
Associated to the contact form is the {\it contact structure}, which is the hyperplane
distribution $\xi\subset TW$ defined to be the kernel of $\lambda$.
We denote such a contact manifold by $(W, \lambda)$. There are three concepts associated to $(W, \lambda)$:
\begin{itemize}
\item The {\it symplectization} of $W$ is, by definition, $W\times\mathbb R$ with symplectic form
$\omega=d(e^t\lambda)$, where $t$ is the coordinate of the factor $\mathbb R$.
We say an almost complex structure $J_\infty$ on $W\times\R$,
is {\it admissible} if it satisfies
\begin{equation}\label{J-infty-2}
\left\{
\begin{array}{ccc}
J_\infty|_\xi & = & J_0,\\
J_\infty\frac{\partial}{\partial t}& =& R_\lambda
\end{array}
\right.
\end{equation}
on $W\times\R$, where $J_0$ is any compatible complex structure on
the symplectic bundle $(\xi,d\lambda)$, $R_\lambda$ is the Reeb vector 
field associated to $\lambda$ defines in (\ref{Reebfield}) below.
Denote by $\mathcal {J}(\lambda)$ the
set of admissible almost complex structures on $W\times\R$.

\item {\it Symplectic completion.} The concept of symplectization
can be generalized to the case of symplectic manifolds with contact type boundary.
Suppose $M$ is a symplectic manifold with contact type boundary $W=\partial M$.
$M$ has a {\it symplectic completion}, which is $$M\cup_{id:W\to W\times\{0\}}(W\times\mathbb R^{\ge 0}),$$
and is denoted by $\widehat{M}$. If $M$ is an exact symplectic manifold, then
$\widehat{M}$ is also exact whose symplectic form
$\widehat{\omega}$ is induced from $M$ and $W\times\mathbb R^{\ge 0}$. In precise,
\begin{equation}\label{cpl-form}
\widehat{\omega}:=
\left\{
\begin{array}{ll}
\omega, & {\rm on}\ M,\\
d(e^t\lambda),& {\rm on}\ W\times\R^+.
\end{array}
\right.
\end{equation}
$M$ is also called the {\it symplectic filling} of $W$ or $W\times\mathbb R$.
Let $J$ be a time-independent almost complex structure on
$\widehat{M}$ which is compatible with $\widehat{\omega}$ and whose
restriction $J_\infty=J|_{W\times\R^+}$ is in $\mathcal
{J}(\lambda)$ and is translation invariant. We denote the space of such $J$ by $\mathcal
{J}(\lambda,\widehat{\omega})$. By \cite{BEHWZ}, such $(\widehat{M},J)$ is an
{\it almost complex manifold with symmetric cylindrical ends adjusted} to the symplectic form $\widehat{\omega}$.

\item A closed {\it Reeb orbit} in $W$ is a closed orbit of the Reeb vector field $Y$:
\begin{equation}\label{Reebfield}
\lambda(Y)=1,\quad\iota(Y)d\lambda=0.
\end{equation}
If the contact form $\lambda$ is
generic, then the set of closed Reeb orbits is discrete.
A closed Reeb orbit $\gamma$ is {\it transversally nondegenerate} if
$$
\det(\mathbb I-d\phi_\lambda^T(\gamma(0))|_\xi)\ne 0,
$$
where $T$ is the period of $\gamma$, $\phi^T_{\lambda}$ is the time-$T$ map of the Reeb
flow. If $\lambda$ is generic, we may
assume all closed Reeb orbits are transversally nondegenerate in $W$.
\end{itemize}

\begin{definition}[Conley-Zehnder Index]
It is easy to see $(\xi, d\lambda|_\xi)$ is a symplectic bundle over $M$, and henceforth,
to each transversally nondegenerate closed Reeb
orbit $\gamma$, one may assign
the corresponding Maslov index, called the {\it Conley-Zehnder index} and
denoted by $\mu_{CZ}(\gamma)$ or simply $\mu(\gamma)$ (see \cite{EGH}).
\end{definition}

A $k$-th iterate $\gamma^k$ of a simple closed Reeb orbit $\gamma$ is {\it good} if
$\mu_{CZ}(\gamma^k)\equiv\mu_{CZ}(\gamma)\mod 2$.
Denote the set of transversally nondegenerate good Reeb orbits by $\mathcal{P}_\lambda$.
We regard a closed Reeb orbit and a multiple of it as two different orbits. And for a closed
orbit $\gamma$, denote by $\kappa_\gamma$ its multiplicity, and
assign the grading $|\gamma|=\mu_{CZ}(\gamma)+n-3$.

\subsection{Pseudo-holomorphic curves}

Let $(\widehat{S},j)$ be a compact smooth oriented surface with a fixed conformal structure $j$, and
$\Lambda=\Lambda^-\sqcup\Lambda^+$ a finite set of {\it interior puncture points}.
We call $S=\widehat{S}\setminus \Lambda$ a punctured Riemann
surface, and the points of $\Lambda^-$ (resp.  $\Lambda^+$) its
incoming (resp. outgoing) points at infinity.

Suppose $W$ is a contact manifold.
There is a deep relationship between the Reeb orbits in $W$
and pseudo-holomorphic curves in $W\times\mathbb R$.
Namely, suppose $u:S\to W\times\mathbb R$ is a pseudo-holomorphic curve.
A theorem of Hofer (see \cite{Hofer} as well as \cite{HWZ}) says that if $u$ is of {\it finite energy}
and has non-removable singular points (punctures),
then these singular points can only approach the Reeb orbits in $W$ at $\pm\infty$.
Let us explain in more detail. Denote by ${\mathcal C}^\pm=\R^\pm\times S^1$ the semi-infinite cylinders.

\begin{definition}\label{cyl-end}
A set of {\it cylindrical ends} for $S$ consists of proper holomorphic embeddings 
$\epsilon_\eta:{\mathcal C}^\pm\to S$, one for each $\eta\in \Lambda^\pm$, using  locally complex coordinates on $S$, 
satisfying
\begin{eqnarray}
\ \epsilon_\eta(r,\theta)=e^{\mp(r+{\rm i}\theta)}\ \ \ {\rm and}\ \ \ \lim_{r\to\pm\infty}\epsilon_\eta(r,\cdot)=\eta
\end{eqnarray}
and with the additional requirement that the images of the $\epsilon_\eta$ are pairwise disjoint.
\end{definition}

\begin{definition}\label{def-asymp}
(1) Suppose $D\subset \mathbb C$ is the unit disk. We say a smooth map
$F:D\setminus \{0\}\to W\times\R$ is {\it asymptotic} to a Reeb orbit $\gamma\subset W$ at $\pm\infty$ if
$ F(r,\theta)=(f(r,\theta),a(r,\theta))$ has the property that
$\lim\limits_{r\to 0}a(r,\theta)=\pm\infty$
and the uniform limit $\displaystyle \lim_{r\to 0}f(r,\theta)$ exists and parameterizes $\gamma$.

(2) More generally, suppose $\Lambda=\{\eta_1,\eta_2,\cdots,\eta_k\}$. 
We say that a smooth map $F:S\to W\times\mathbb R$ is {\it asymptotic}
to a Reeb orbit ${\gamma}_j\in \mathcal{P}_{\lambda}$ in $\eta_j$ at $\pm\infty$, if there exists polar
coordinates $(r,\theta)$ centered at $\eta_j$ such that $F$ restricted to a neighborhood of $\eta_j$
is asymptotic to $\gamma_j$ in the sense above. These $\eta_j$'s are called the {\it punctures} of $F$;
they are either positive or negative according to the sign of $\pm\infty$.
\end{definition}

Suppose $F=(f,a):S\to W\times \mathbb R$ is a pseudo-holomorphic curve, {\it i.e.} $dF\circ j=J_{\infty}\circ dF$,
with the set of punctures, say $Z$. The {\it energy} of $F$ is defined to be
$$
E_{\lambda}(F)=\sup_{\varphi\in\mathcal C}\int_{S\backslash Z} (\varphi\circ a)da\wedge f^*\lambda,
$$
where $\mathcal C$ is the set of all non-negative smooth functions $\varphi:\mathbb R\to\mathbb R$
having compact support and satisfying the condition $\int_{\mathbb R}\varphi(x)dx=1$.

\begin{theorem}[Hofer \cite{Hofer}]
Suppose $Z=\{\eta_1,\eta_2,\cdots,\eta_k\}\subset S$
is a finite subset of a Riemann surface $S$. Then every
pseudo-holomorphic curve $F:S\backslash Z\to W\times \mathbb R$
of finite energy and without removable singularities is asymptotic
to a closed Reeb orbit $\gamma_i$ in $W$ near each puncture $\eta_i$.
\end{theorem}

\subsection{The moduli spaces and their compactification and orientation}\label{orient_comp}

Choose a generic contact form such that the Reeb orbits are discrete.
With Hofer's theorem, one may consider
the moduli space of pseudo-holomorphic curves with a given type. Namely, let $\widehat{S}$ be a Riemann surface,
$S=\widehat{S}\setminus \{\gamma_1^+,\cdots,\gamma_s^+\}\cup\{\gamma_1^-,\cdots,\gamma_t^-\}$,
and let
$\widehat{\mathcal M}(S; \gamma_1^+,\cdots,\gamma_s^+;\gamma_1^-,\cdots,\gamma_t^-)$ 
be the set of punctured
pseudo-holomorphic curves
\begin{equation}\label{fun_pseudo}
F=(f,a):S\to W\times\mathbb R
\end{equation}
such that near the
positive punctures $q_i^+$ and negative punctures $q_j^-$, the pseudo-holomorphic curve
is asymptotic
to Reeb orbits $\gamma^+_i$ and $\gamma^-_j$, respectively.

The compactification of the moduli spaces in SFT is studied in \cite{BEHWZ}. Roughly, denote by
$${\mathcal M}(S; \gamma_1^+,\cdots,\gamma_s^+;\gamma_1^-,\cdots,\gamma_t^-)=
\widehat{\mathcal M}(S; \gamma_1^+,\cdots,\gamma_s^+;\gamma_1^-,\cdots,\gamma_t^-)/ \mathrm{Aut}(S)$$
the moduli space of punctured pseudo-holomorphic curves.
There is a natural $\mathbb R$ action on $\mathcal M$ which shifts the pseudo-holomorphic
curves vertically. By modulo such an $\mathbb R$ action, it
is proved in symplectic field theory (\cite{BEHWZ}) that
the space of pseudo-holomorphic curves with a given type and of finite energy (the bound of the energy is
{\it a priori} chosen)
can be partially compactified into a stratified space, whose codimension greater than zero
strata is described by the ``broken'' curves, which, topologically,
can be realized as pseudo-holomorphic curves which are stretched
to be infinitely long in the middle of $W\times \R$ at some time.

The orientation issue is fully discussed in \cite{BM}. Similar to the Hamiltonian chord case,
to each closed Reeb orbit $\gamma$,
there is an associated
Fredholm operator (see \cite[\S2 Proposition 4]{BM}).
Its index is exactly $|\gamma|$ and its oriented bundle $o_{\gamma}$
is the determinant bundle of this operator. The orientation bundle of
${\mathcal M}(S; \gamma_1^+,\cdots,\gamma_s^+;\gamma_1^-,\cdots,\gamma_t^-)$
is then
$$\det S\otimes
o_{\gamma_1^+}^-\otimes\cdots\otimes o_{\gamma_s^+}^-
\otimes o_{\gamma_1^-}^-\otimes\cdots\otimes
o_{\gamma_t^-}^-.$$
Note that, if we switch $\gamma_i^{+}$ with $\gamma_{i+1}^{+}$,
then the orientation of
${\mathcal M}(S; \gamma_1^+,\cdots,\gamma_s^+;\gamma_1^-,\cdots,\gamma_t^-)$
and the one of
${\mathcal M}(S; \gamma_1^+,\cdots,\gamma_{i+1}^+,\gamma_i^+,\cdots,
\gamma_s^+;\gamma_1^-,\cdots,\gamma_t^-)$ agree if and only
if $|\gamma_i^{+}||\gamma_{i+1}^+|$ is even.
For more details, see \cite{BM}.

\subsection{Linearized contact homology}
The theory of pseudo-holomorphic curves in a symplectization $W\times\mathbb R$
can be generalized to the case of a symplectic completion
$\widehat M=M\cup_{id:W\to W\times\{0\}}W\times\mathbb R^{\ge 0}$
for a symplectic manifold $M$ with contact type boundary.
The good property to consider this case is that, the theory of pseudo-holomorphic
curves in $W\times\mathbb R$ gives rise to the theory of ``contact homology'' of $W$ (see \cite{EGH}),
while if we consider both, the theory of contact homology is richer and admits a ``linearization''.

\begin{definition}[Augmentation]
Given $\gamma\in\mathcal{P}_\lambda$,
denote by $\mathcal {M}(\gamma;\emptyset )$ the moduli space of (equivalence classes of)
$J$-holomorphic planes $F:\mathbb C\to\widehat{M}$
which is asymptotic to $\gamma$ at $\infty$. If $J$ is regular and
$\dim\mathcal{M}(\gamma;\emptyset )=0$, denote by
$$\varepsilon(\gamma):=\#\mathcal{M}(\gamma;\emptyset ),
\quad\mbox{for all }\gamma\in\mathcal{P}_\lambda,$$
where $\#$ is counting with sign (see \cite{BM}),
and $\varepsilon(\gamma)$ is called the {\it symplectic augmentation}
of $\gamma$.
\end{definition}

Let $\Cont_*^{\mathrm{lin}}(M)$ be the linear space spanned by
$\mathcal{P}_\lambda$ over a field $\mathbb K$ of
characteristic zero,
graded as above. Under our exactness conditions, we can  define a linear operator
$$d:\Cont_*^{\mathrm{lin}}(M)\to \Cont_{*-1}^{\mathrm{lin}}(M)$$
by
\begin{equation}\label{def_contdiff}
d(\gamma^+):=\sum_{
\genfrac{}{}{0pt}{}{\gamma^-, \gamma_1^-,\cdots,\gamma_t^-:}{|\gamma^-|+|\gamma_1^-|+\cdots+|\gamma_t^-|=|\gamma^+|-1}
} \frac{\#(\mathcal M (\gamma^+; \gamma^-,\gamma_1^-,\cdots,\gamma_t^-)/\mathbb R)}
{\kappa_{\gamma^-}\cdot\kappa_{\gamma_1^-}\cdots\kappa_{\gamma_t^-}}
\varepsilon(\gamma_1^-)\cdots\varepsilon(\gamma_t^-)\cdot \gamma^-,
\end{equation}
where $\mathcal M (\gamma^+; \gamma^-,\gamma_1^-,\cdots,\gamma_t^-)$
is the moduli space of pseudo-holomorphic {\it spheres} with punctures asymptotic
to $\gamma^+$ at $+\infty$ and $\gamma^-,\gamma_1^-,\cdots,\gamma_t^-$ at $-\infty$,
respectively. 

\begin{definition-lemma}[Linearized Contact Homology]
Let $M$ be a symplectic manifold with contact boudary, and
let $d:\Cont_*^{\mathrm{lin}}(M)\to \Cont_{*-1}^{\mathrm{lin}}(M)$ be defined by $(\ref{def_contdiff})$.
Then $d^2=0$, and the associated homology is
called the {\it linearized contact homology} of $M$, denoted by $\mathrm{CH}^{\mathrm{lin}}_*(M)$.
\end{definition-lemma}

Note that the general definition of $d$ is a refined version of (\ref{def_contdiff}) involving the 
Novikov ring, and the proof heavily depends on the polyfold theory which is currently being
developed by Hofer and his collaborators. We shall not discuss it here, and
the interested reader may refer to Cieliebak-Latschev \cite[\S5]{CL} for an algebraic exposition,
and also \cite[\S3]{BO} for more details.

\subsection{Lie bialgebra of Cieliebak-Latschev}

Cieliebak-Latschev proved in \cite{CL} that the linearized contact homology
of $M$ endows the structure of an involutive Lie bialgebra. More precisely, they proved
that {\it on the chain level}, the linearized contact chain complex $\mathrm{Cont}_*^{\mathrm{lin}}(M)$
forms what they called a ``BV$_\infty$ algebra''.

A BV$_\infty$ algebra is in many aspects similar to a bi-Lie$_\infty$ algebra.
It contains a (graded) Lie$_\infty$ algebra and a (graded) Lie$_\infty$ coalgebra, with some
compatibility conditions. In general, a Lie bialgebra
may not necessarily be involutive, and similarly, on the chain level,
a bi-Lie$_\infty$ algebra may not be involutive even ``up to homotopy''.
However, a BV$_\infty$ algebra is {\it a priori} involutive up to homotopy.
The homotopy for involutivity has a deep relation with the structure on the
moduli space of Riemann surfaces of all genera, which we prefer not to discuss in current paper.

While the whole theory of BV$_\infty$ is to be developed by Cieliebak-Latschev
(see, however, \cite{CL}), in the following we briefly introduce part of their results, 
the ones that are simple algebraically and geometrically.

First, we introduce the definition of Lie$_\infty$ algebras.
More details on this concept can be found at
Lada-Stasheff \cite{LS}.

Suppose $L$ is a graded vector space over $\mathbb K$. Let $\ov L$
be the desuspension of $L$, and let $\bigwedge^\bullet L$ be the
graded symmetric tensor algebra generated by $\ov L$, which may be identified with the exterior algebra
generated by $L$. We would like to view $\bigwedge^\bullet L$
as a cocommutative coalgebra instead of a commutative algebra, where the coproduct is given by
$$
\Delta(a_1 a_2\cdots  a_n)=\sum_{p+q=n}\sum_{\sigma}\pm
a_{\sigma(1)}a_{\sigma(2)}\cdots a_{\sigma(p)}
\otimes a_{\sigma(p+1)} a_{\sigma(p+2)}\cdots a_{\sigma(p+q)},
$$
where $\sigma$ runs over all $(p,q)$-unshuffles of $n$.

\begin{definition}[Lie$_\infty$ Algebra]\label{def_lieinfty}
Suppose $L$ is a graded vector space over $\mathbb K$.
A Lie$_\infty$ algebra on $L$ is a degree $-1$ differential
$\delta:\bigwedge^\bullet L\to\bigwedge^\bullet L$
which is a coderivation with respect to the coproduct $\Delta$.
\end{definition}

The coderivation $\delta$ can be expanded by its ``Taylor series''.
More precisely, a Lie$_\infty$ algebra consists of a sequence
of linear operators
$$\delta_n:\mbox{$\bigwedge^n L$}\to\ov L,\quad n=1,2,\cdots$$
such that
\begin{equation}\label{lieinfty_reln}
\sum_{i+j=n}\sum_{\sigma}\pm \delta_{j+1}(\delta_i(a_{\sigma(1)}\cdots
a_{\sigma(i)})a_{\sigma(i+1)}\cdots a_{\sigma(n)})=0,
\end{equation}
for all $a_1\cdots a_n\in\bigwedge^n L$, where $\sigma$ runs over all
$(i, j)$-unshuffles of $n$.
In equation (\ref{lieinfty_reln}), if we
apply $\delta_n$ to $\bigwedge^m L$ by
$$
\delta_n(a_1a_2\cdots a_m)=\left\{\begin{array}{cl}\displaystyle\sum_{\sigma}\pm
\delta_n(a_{\sigma(1)}\cdots a_{\sigma(n)})a_{\sigma(n+1)}\cdots a_{\sigma(m)},& m\ge n,\\
0,&m<n,
\end{array}\right.
$$
and set $\delta:=\delta_1+\delta_2+\cdots$,
then $\delta$ exactly gives a Lie$_\infty$ algebra on $L$.

\begin{example}[Lie Algebra]
A Lie algebra is naturally a Lie$_\infty$ algebra. In fact, suppose $L$ is a Lie algebra over $k$.
The Eilenberg-Chevalley complex of $L$ is $\bigwedge^\bullet L$, with the
differential $\partial$ defined by
$$\partial (a_1a_2\cdots a_n):=\sum_{i<j}\pm [a_i,a_j]a_1\cdots\hat a_i\cdots\hat a_j\cdots a_n.$$
The Jacobi identity implies $\partial^2=0$. If we set $\delta_2=\partial$ and $\delta_i=0$ for all $i\ne 2$,
then the Eilenberg-Chevalley complex of $L$ exactly gives a Lie$_\infty$ algebra structure on $L$.
\end{example}

On the other hand, any Lie$_\infty$ algebra $(L,\delta)$ gives rise to
a Lie algebra on $L$ ``up to homotopy''. Namely, for any $a_1,a_2\in L$,
let $\ov a_1,\ov a_2$ be their image under the desuspension $L\to\ov L$
and set
$$
[a_1,a_2]:=\mbox{suspension of}~(-1)^{|a_1|}\delta_2(\ov a_1,\ov a_2),
$$
then $[\,,\,]$ thus defined is graded skew-symmetric, and
$\delta_1\circ\delta_3+\delta_2\circ\delta_2+\delta_3\circ\delta_1=0$
implies Jacobi identity up to homotopy. That is, $\HH_*(L,\delta_1)$
is a graded Lie algebra.

\begin{theorem}[Cieliebak-Latschev \cite{CL}]
Let $M$ be a symplectic manifold with contact type boundary such
that $c_1(M)=0$. Then
the linearized contact homology $\mathrm{CH}_*^{\mathrm{lin}}(M)$
has the structure of an involutive Lie bialgebra of degree $(2-n,2-n)$. More precisely,
the linearized contact chain complex $\mathrm{Cont}_*^{\mathrm{lin}}(M)$
has the structure of a BV$_\infty$ algebra; in particular, $\Cont_*^{\mathrm{lin}}(M)$
forms a Lie$_\infty$ algebra.
\end{theorem}

\begin{proof}[Sketch of proof]
Denote by $\mathcal{M}(\gamma_1^+,\gamma_2^+;\gamma^-,\gamma_1^-,\cdots,\gamma_t^-)$
 the moduli space of pseudo-holomorphic spheres with two punctures at
$+\infty$ and $t+1$ punctures at $-\infty$ in the symplectization
$W\times\R$.
As in the definition of linearized contact homology,
we want to remove those pseudo-holomorphic curves that are
``$t$-to-$0$''. Note that in the definition of linearized contact chain complex,
since there is only one incoming Reeb orbits,
the only possibility is $1$-to-$0$. In the general case, this might not be true any more.
Similar to symplectic augmentation,
we define
$\varepsilon(\gamma_1,\gamma_2):=\#\mathcal M(\gamma_1,\gamma_2;\emptyset)$
be the number of pseudo-holomorphic spheres in $\widehat M$ with two punctures at $+\infty$.

Let
\begin{eqnarray}\label{def_bracket}
[\gamma_1^+,\gamma_2^+]&:=&\sum_{\gamma^-,\gamma_1^-,\cdots,\gamma_t^-}
\frac{\#(\mathcal M (\gamma_1^+,\gamma_2^+; \gamma^-,\gamma_1^-,\cdots,\gamma_t^-
)/\mathbb R)}{\kappa_{\gamma^-}\cdot\kappa_{\gamma_1^-}\cdots\kappa_{\gamma_t^-}}
\varepsilon(\gamma_1^-)\cdots\varepsilon(\gamma_t^-)\cdot \gamma^-\nonumber \\
&+&\sum_{\gamma^-,\gamma_1^-,\cdots,\gamma_t^-}
\frac{\#(\mathcal M (\gamma_1^+; \gamma^-,\gamma_1^-,\cdots,\gamma_t^-
)/\mathbb R)}{\kappa_{\gamma^-}\cdot\kappa_{\gamma_1^-}\cdots\kappa_{\gamma_{t}^-}}\varepsilon(\gamma_1^-)\cdots
\varepsilon(\gamma_{t-1}^-)\varepsilon(\gamma_2^+,\gamma_{t}^-)\cdot \gamma^-\nonumber\\
&+&\sum_{\gamma^-,\gamma_1^-,\cdots,\gamma_t^-}
\frac{\#(\mathcal M (\gamma_2^+; \gamma^-,\gamma_1^-,\cdots,\gamma_t^-
)/\mathbb R)}{\kappa_{\gamma^-}\cdot\kappa_{\gamma_1^-}\cdots\kappa_{\gamma_t^-}}\varepsilon(\gamma_1^-)\cdots
\varepsilon(\gamma_{t-1}^-)\varepsilon(\gamma_1^+,\gamma_{t}^-)\cdot \gamma^-,
\end{eqnarray}
which counts all pseudo-holomorphic pants in $W\times \mathbb R$
with two sleeves asymptotic to $\gamma_1^+,\gamma_2^+$
at $+\infty$ and one sleeve asymptotic to a closed Reeb orbit $\gamma^-$ at $-\infty$.
The orientation of the moduli spaces (\S\ref{orient_comp})
guarantees that $[\,,\,]$ thus defined is graded skew-symmetric.

To show that $[\,,\,]$ respects the contact differential, one considers the compactification
of the moduli space of the above pseudo-holomophic pants; if it is one dimensional,
then its boundary exactly gives
$$d[\gamma_1^+,\gamma_2^+]
=[d\gamma_1^+,\gamma_2^+]+(-1)^{|\gamma_1^+|}[\gamma_1^+,d\gamma_2^+].$$
Similarly, to show the Jacobi identity holds up to homotopy,
we consider the moduli space of pseudo-holomorphic spheres with three punctures at $+\infty$
and one punctures at $-\infty$; if it is one dimensional, then the boundary of its compactification
gives the Jacobi identity up to homotopy. The homotopy for homotopies of the Jacobi identity, and all
higher homotopies, are given by pseudo-holomorphic curves with all possible punctures at $+\infty$.

The Lie co-bracket is defined similarly:
\begin{equation}\label{def_cobracket}
\delta(\gamma^+):=\sum_{\gamma^-,\tilde\gamma^-,\gamma_1^-,\cdots,\gamma_n^-}
\frac{\#(\mathcal M (\gamma^+; \gamma^-,\tilde\gamma^-,\gamma_1^-,\cdots,\gamma_n^-
)/\mathbb R)}{\kappa_{\gamma^-}\cdot\kappa_{\tilde\gamma^-}\cdot\kappa_{\gamma_1^-}\cdots\kappa_{\gamma_n^-}}
\varepsilon(\gamma_1^-)\cdots\varepsilon(\gamma_n^-)\cdot \gamma^-\wedge\tilde\gamma^-.
\end{equation}
The proof that $[\,,\,]$ and $\delta(\,)$ form a Lie bialgebra on the homology level is given in \cite{CL},
where the readers may find more interesting structures.

Finally, a word of degrees: the linearized contact homology comes
from the {\it linearization} of the contact homology, originally defined
in \cite{EGH}. The chain complex for the contact homology, is a free DG algebra generated
by the closed Reeb orbits, together with
a degree $-1$ differential, which is giving by counting pseudo-homomorphic spheres with
punctures at $\pm\infty$. If one re-grade the closed Reeb orbits by their Conley-Zehnder index
shifted by one, then we do get that the degree of the Lie bialgebra is $(2-n,2-n)$.
More details can be found in \cite{CL}.
\end{proof}

\section{Analytic setting of pseudo-holomorphic disks with punctures}\label{sect_ana}

In the rest of the paper we are going to define the chain map
$f:{\rm Cont}_*^{\mathrm{lin}}(M)\to{\rm Cycl}^*(\mathcal Fuk(M))$. In order
to do that, we need to
study the punctured pseudo-holomorphic curves with boundaries in
Lagrangian submanifolds, and near each puncture the curves are asymptotic to some periodic Reeb orbit.
Moduli spaces that appear in Fukaya category and in symplectic field theory
have been well studied ({\it c.f.} \cite{FOOO,Seidel} and \cite{BEHWZ,BM}),
however, moduli spaces that combine these two are less studied in literature,
and for the sake of completeness, we discuss them in a little more detail.
However, both the compactifications and the orientations of the
moduli spaces of these curves are the combinations of the 
compactifications and orientations of those studied in SFT and Fukaya category respectively.

\subsection{Basic setting}

We fix a field $k$. Before we deal with the orientability issue, we only restrict to the
case of $\mbox{Char}(k)=2$. Recall that in our specific case, we only
consider an exact symplectic manifold
$M$ with contact type boundary $W$ such that $c_1(M)=0$, 
and its admissible Lagrangian submanifolds
(see Assumption \ref{assum-Fuk} in \S\ref{fukcat}). 
Also, we only consider transversally nondegenerate {\it good} Reeb periodic orbits.
Let us first recall some definitions and notations.

\subsubsection{Punctured pointed-boundary Riemann surfaces}

Let $(\widehat{S},j)$ be a compact smooth oriented surface with boundary
and with a fixed conformal structure $j$, and $\Sigma$
a finite set of {\it boundary points}, divided into two parts $\Sigma=\Sigma^-\sqcup\Sigma^+$, and
$\Lambda=\Lambda^-\sqcup\Lambda^+$ a finite set of {\it interior puncture points}.
We call $S=\widehat{S}\setminus (\Sigma\sqcup\Lambda)$ a punctured pointed-boundary Riemann
surface, and the points of $\Sigma^-$ and $\Lambda^-$ (resp. $\Sigma^+$ and $\Lambda^+$) its
incoming (resp. outgoing) points at infinity. For instance, we use the following special notations for
some simpler surfaces: (1) $D$ for the closed unit disc in $\mathbb C$; (2) $H$ (resp. $\ov{H}$) for the
closed upper half plane, with one incoming (resp. outgoing) point at infinity; (3) $Z=\R\times [0,1]$
for infinite strip with the coordinates $(s,t)$; (4) $\mathcal C=\R\times S^1$ for infinite cylinder
with the cylinder coordinates $(r,\theta)$. Both $Z$ and $\mathcal C$
have an incoming point $s,r=-\infty$ and an outgoing point $s,r=+\infty$.

Recall the definition of Lagrangian labels (see \S\ref{PbRS}):
A set of Lagrangian labels for $S$ is a family $\mathcal{L}:=\{L_C\}$ of 
admissible Lagrangian submanifolds $L_C\subset M$, indexed by the 
connected components $C\subset \partial S$. Each 
$\zeta\in\Sigma\subset \widehat{S}$ is in the closure of two boundary 
components, say $C_{\zeta,0}$ and $C_{\zeta,1}$. 
When considering maps 
from $S$ to $M$, $C_{\zeta,0}$ and $C_{\zeta,1}$ are mapped 
into $L_{C_{\zeta,0}}$ and $L_{C_{\zeta,1}}$ respectively.
If $M$ has a contact type boundary, then each 
such $L_C$ is a Lagrangian submanifold in $(\widehat{M},\widehat{\omega})$, 
and to distinguish these two sets of Lagrangian labels, we denote
the former by $(M,\mathcal L)$ and later by $(\widehat{M},\mathcal{L})$.



\subsubsection{Punctured pointed-disks}

We first study the moduli space of pseudo-holomorphic disks with one interior puncture. 
For cases of several interior punctures, it will be clear that only minor modifications are needed. 
A $(1,k)$-disk $S$ is a punctured
pointed-boundary Riemann surfaces whose compactification $\widehat{S}\cong D^2$, and with 1 negative (incoming)
interior puncture $\eta$ and $k$ positive (outgoing) boundary punctures $\zeta_1,\cdots,\zeta_k$.
Number the marked points on the boundary respecting their (anti-clock wise) cyclic order along the boundary
(induced by the orientation on the disk), and denote the corresponding strip-like ends by
$\varepsilon_1,\cdots,\varepsilon_k$ and cylindrical end still by $\varepsilon_\eta$.
Denote by ${\mathcal R}^1_k$ the moduli space of such disks.
The moduli space of $(k+1)$-pointed-boundary disks (without interior puncture), equipped
with 1 incoming strip-like end $\epsilon^-_0$ and $k$ outgoing strip-like ends
$\epsilon^+_1,\cdots,\epsilon^+_k$, is denoted by
${\mathcal R}_k.$
The codimension one strata of the Deligne-Mumford compactification
$\overline{\mathcal R}^1_k$ are $\bigcup_{k_1+k_2=k+1}
{\mathcal R}^1_{k_1}\times{\mathcal R}_{k_2}$, $2\le k_2\le k$.

\subsubsection{Gluing of domains}\label{glue_of_domains}

(1) We can explicitly describe the {\it gluing} process of a 1-punctured
pointed-boundary surface $S_1\in\mathcal{R}^1_{k_1}$ and a
pointed-boundary surface $S_2\in\mathcal{R}_{k_2}$ at a
boundary marked point $\zeta_1\in\Sigma_1^+$ with a positive
strip-like end $\epsilon^+_{\zeta_1}$ and $\zeta_2\in\Sigma_2^-$
with a negative strip-like end $\epsilon^-_{\zeta_2}$, as follows.
For any {\it gluing length} $l>0$, let
$S_1^*=S_1\setminus \epsilon^+_{\zeta_1}((l,+\infty)\times[0,1])$,
$S_2^*=S_2\setminus \epsilon^-_{\zeta_2}((-\infty,-l)\times[0,1])$,
if we identify $\epsilon^+_{\zeta_1}(s,t)\sim\epsilon^-_{\zeta_2}(s-l,t)$
for $(s,t)\in[0,l]\times[0,1]$, then we obtain an $l$-length glued surface
$S=S_1\#_l S_2=(S^*_1\cup S^*_2)/\sim$, which is in
$\mathcal{R}^1_{k_1+k_2-1}$. Gluing the two conformal structures
$j_1$ on $S_1$ and $j_2$ on $S_2$, we can construct glued
conformal structure $j_l$ on  $S_1\#_l S_2$.

(2) Similarly, we can glue  another punctured surface (with or without
boundary and with at least one positive cylindrical end) and a
1-punctured pointed-boundary surface at the interior punctures.
As a special example, we
consider the gluing process of a 2-punctured sphere
$\tilde{S}^2=S^2\setminus\{\eta^+,\eta^-\}$ with punctures
$\eta^+$ and $\eta^-$ (one with positive cylindrical end and
the other with negative cylindrical end)  and a 1-punctured
pointed-boundary surface $S\in\mathcal{R}^1_k$ at the interior
puncture $\eta$ with a negative cylindrical end $\epsilon^-_{\eta}$.
First, the 2-punctured sphere $\tilde{S}^2$ can be conformally
equivalent to a cylinder $\tilde{S}^2=h(Z)$ such that $\lim_{r\to\pm\infty}h(r,\cdot)=\eta^\pm$.
Then take gluing length $l>0$,
let $(\tilde{S}^2)^*=\tilde{S}^2\setminus h((l,+\infty)\times S^1)$,
$S^*=S\setminus \epsilon^-_{\eta}((-\infty,-l)\times S^1$, then
identify $h(r,\theta)\sim \epsilon^-_{\eta}(r-l,\theta)$ for
 $(r,\theta)\in[0,l]\times S^1$, thus we obtain a $l$-length glued
 surface $S'=\tilde{S}^2\#_l S=((\tilde{S}^2)^*\cup S^*)/\sim$,
 which is still in $\mathcal{R}^1_k$. Similarly, we can construct
 glued conformal structure $j_l$ on $\tilde{S}^2\#_l S$.
The general case is treated similarly.

\subsubsection{Analytic data}

Recall $\mathcal
{J}(\lambda,\widehat{\omega})$ is the space of time-independent almost complex structure $J$ on
$\widehat{M}$ which is compatible with $\widehat{\omega}$ and whose
restriction $J_\infty=J|_{W\times\R^+}$ is in $\mathcal
{J}(\lambda)$ and is translation invariant.

\begin{definition}[Analytic Data for Relating Maps]\label{def-data-relat}
Let $S\in\overline{\mathcal{R}}^1_k$ be a stable one-punctured pointed-boundary
Riemann surface with  (compatibly labeled) Lagrangian labels. The
{\it analytic data for relating maps} on $S$, denoted by $\mathbf D_{rel}$, consists of the following choices:

\begin{itemize}
\item Cylindrical end $\epsilon^-_\eta$ for incoming interior puncture and strip-like
ends $\epsilon^+_1,\cdots,\epsilon^+_k$ for outgoing boundary punctures.

\item  A Floer datum $(H_\zeta,J_\zeta)$ for each pair of Lagrangian submanifolds
$(L_{\zeta,0},L_{\zeta,1})$ associated to  $\zeta\in\Sigma$ is as above;
\item A perturbation datum for $S$ (with or without punctures) is a pair $(K,J)$, where
$K\in\Omega^1(S,\mathscr H)$ is a function-valued 1-form on $M$ satisfying
$K(\xi)|_{L_C}=0$ for all $\xi\in TC\subset T(\partial S)$; $J\in C^\infty(S,\mathscr J)$
is a family of almost complex structure; moreover, $K$ and $J$ should be {\it compatible}
with the chosen cylindrical and strip-like ends and Floer data, $i.e.$
\begin{eqnarray}
\epsilon^*_\eta K=0,&& J(\epsilon_\eta(r,\theta))=J_\eta(\theta),\label{comp-cond-1}\\
\epsilon^*_{\zeta}K=H_\zeta(t)dt, && J(\epsilon_\zeta(s,t))=J_\zeta(t)\label{comp-cond-2}
\end{eqnarray}
for each $\eta\in\Lambda^\pm$, $(r,\theta)\in \mathcal{C}^\pm$, and $\zeta\in\Sigma^\pm$, $(s,t)\in Z^\pm$.
\end{itemize}
\end{definition}

Note that $K$ determines a vector-field-valued 1-form $Y\in\Omega^1(S,C^\infty(TM))$:
for each $\xi\in TS$, $Y(\xi)$ is the Hamiltonian vector field of $K(\xi)$.

\begin{example} We can identify $Z=\{x+iy\in\mathbb C: x\in\mathbb R, y\in[0,1]\}$ which is
conformal  to $D^2\backslash \{\pm 1\}\subset\mathbb C$. Label the upper boundary $Z^{+}$
of $Z$ by $L_0$, and the lower boundary $Z^{-}$ by $L_1$, where $L_0, L_1$ are two admissible
Lagrangian submanifolds. A  perturbation datum for $Z$ (without puncture) is a pair $(K, J)$,
where $K\in\Omega^1(Z, \mathscr H)$ satisfies
$K(\xi)|_{L_{1,2}}=0, \mbox{ for all }\xi\in TZ^{+} \cup TZ^{-},$
and $J$ is a $Z$-family of complex structures in $\mathscr J$. In addition,
$K$ and $J$ are compatible with the Floer data as given in (\ref{comp-cond-2}).
\end{example}

Since the gluing operation involves pointed-disks may {\it a priori} produce different sets of strip-like
ends or/and  perturbation data, we have to get a careful choice as follows. Note that we have the {\it gluing map}
$$
\begin{array}{cccl}
\gamma:&\overline{\mathcal{R}}_{\mathrm{glue} }=(0,+\infty]^m\times(\partial\mathcal{R}_k)^m&\longrightarrow&
\overline{\mathcal{R}}_k\\
&\prod_{i=1}^m l_i\times\prod_{j=1}^{m+1}S_j&\longmapsto& S,
\end{array}
$$
where $(\partial\mathcal{R}_k)^m$ is the set of codimension $m$ strata of $\overline{\mathcal{R}}_k$,
$l_i$ are the gluing lengths. $\mathcal{R}_{\mathrm{glue}}=(0,+\infty)^m\times(\partial\mathcal{R}_k)^m$ is
called the {\it glued space}.

Similarly, for the space $\mathcal{R}^1_k$, we have the {\it gluing map}
$$
\begin{array}{cccl}
\beta:&\overline{\mathcal{R}}^1_{\mathrm{glue}}=(0,+\infty]^m\times(\partial\mathcal{R}^1_k)^m&\longrightarrow&
\overline{\mathcal{R}}^1_k\\
&\prod_{i=1}^m l_i\times\prod_{j=1}^{m+1}S_j&\longmapsto&S,
\end{array}
$$
where $(\partial\mathcal{R}^1_k)^m$ is the set of codimension $m$ strata of
Deligne-Mumford compactification $\overline{\mathcal R}^1_k$, for instance,
$\displaystyle (\partial\mathcal{R}^1_k)^1=\bigcup_{k_1+k_2=k+1}
{\mathcal R}^1_{k_1}\times{\mathcal R}_{k_2}$, $2\le k_2\le k$,
$l_i$ are the gluing lengths (see case (1) of \S\ref{glue_of_domains}). $\mathcal{R}^1_{\mathrm{glue}}
:=(0,+\infty)^m\times(\partial\mathcal{R}^1_k)^m$
is called the {\it glued space}. Note that the analytic data on some
boundary stratum might be coming from the analytic data $\mathbf D_{\mathrm{Fuk}}$
from Fukaya category.

\begin{definition}[Consistent Universal Analytic Data for Relating Maps]\label{Consist-relate}
A consistent universal choice of analytic data for a relating map is a choice $\mathbf D_{rel}$ of
analytic data,
for each integer $k\geq 1$ and every  representative of $S\in\overline{\mathcal R}^1_k$,
such that:

(1) For any possible $m$, there is an open subset $\overline{U}\subset\overline{\mathcal{R}}^1_{\mathrm{glue}}$
containing the $\{+\infty\}\times(\partial\mathcal{R}^1_k)^m$, such that the two analytic data
(coming from $\mathbf D_{rel}$ or $\mathbf D_{\mathrm{Fuk}}$) on the glued space (one is inherited from the
universal choice of data on each boundary stratum through the gluing process, the other
is the pullback from the universal choice of data on $\mathcal{R}^1_k$ via the gluing map)
agree over $U=\overline{U}\cap\mathcal{R}^1_{\mathrm{glue}}$;

(2) Let $(K,J)$ be the first perturbation datum (obtained by gluing) on $\mathcal{R}^1_{\mathrm{glue}}$, and
$(\overline{K},\overline{J})$ its extention to the compactification $\overline{\mathcal{R}}^1_{\mathrm{glue}}$,
then the other datum (obtained by pullback  from $\mathcal{R}^1_k$) also extends smoothly
to $\overline{\mathcal{R}}^1_{\mathrm{glue}}$, and the extension agrees with $(\overline{K},\overline{J})$
over the subset $\{+\infty\}\times(\partial\mathcal{R}^1_k)^m\subset\overline{\mathcal{R}}^1_{\mathrm{glue}}$.
\end{definition}

\begin{lemma}
Consistent universal choices of analytic data for both Fukaya category and relating maps exist.
\end{lemma}
\begin{proof}
The argument is the same as the one in Lemmas 9.3 and 9.5 of \cite{Seidel}.
\end{proof}

\subsubsection{Inhomogeneous pseudo-holomorphic maps}

Suppose $S$ is an element in $\mathcal R^1_k$ with Lagrangian labels. Equip it with
cylindrical and strip-like ends, and consistent universal analytic data {\bf D}$_{rel}$.

Consider the {\it inhomogeneous pseudo-holomorphic map} equation for
$u\in C^{\infty}(S, \widehat M)$ which is positive asymptotic to time-1 Hamitonian chord $y_{\zeta_i}$ at each
$\zeta_i$, $i=1,\cdots,k$, and negative asymptotic to Reeb orbit $\gamma$ at puncture (see Figure~4):
\begin{equation}\label{puncture-ipm}
\left\{\begin{array}{l}
Du(z)+J(z,u)\circ Du(z)\circ j_S=Y(z,u)+J(z,u)\circ Y(z,u)\circ j_S,\\
u(C)\subset L_C,\quad {\rm for\ all\ } C\subset\partial S,\\
\displaystyle\lim_{s\to+\infty}u\circ\epsilon_{\zeta_i}(s,\cdot)=
y_{\zeta_i}, \quad {\rm for}\ \zeta_i\in\Sigma^+,\ y_{\zeta_i}\in\Hom(\tilde L_{\zeta_i,0},\tilde L_{\zeta_i,1}),\\
\displaystyle\lim_{r\to-\infty}F_\eta(r,\cdot)=+\infty,\ \ \
\displaystyle\lim_{r\to-\infty}a(r,\cdot)=+\infty,\\
\displaystyle\lim_{r\to-\infty}f(r,\cdot)=\gamma\in\mathcal{P}_\lambda,\ \ \
\displaystyle\lim_{z\to 0,z\in \ell}f(z)=m_\gamma\in\gamma.
\end{array}
\right.
\end{equation}
where $j_S$ is the complex structure on $S$, $J\in\mathscr J$,
$F_\eta(r,\theta)=(f(r,\theta),a(r,\theta))=u\circ\epsilon_\eta(r,\theta)$.
We denote by $\widetilde{\mathcal {M}}(\gamma;\{y_{\zeta_i}\}_{i=1}^k; \widehat{M},\mathcal{L})$ the
set of tuples $(j_S,\eta,\vec{\zeta},u)$, where $S\in\mathcal{R}^1_k$, $u$ is the solution of
(\ref{puncture-ipm}), and $\vec{\zeta}$ denotes the ordered boundary punctures $(\zeta_1,\cdots,\zeta_k)$.
For fixed $S$, and so fixed $(j_S,\eta,\vec{\zeta})$, denote by $\widetilde{\mathcal {M}}_S
(\gamma;(y_{\zeta_1},\cdots,y_{\zeta_k}); \widehat{M},\mathcal{L})$ the set of solution $u$ of (\ref{puncture-ipm}).

Given a 1-punctured $k$-pointed-boundary disc $S$ with complex structure $j_S$.
Denote the group of automorphism on the domain by $\mathrm{Aut}(S)$.
Then the {\it moduli space of 1-punctured $k$-pointed-boundary} ($(1,k)$-punctured, for short) 
{\it pseudo-holomorphic maps} of $S$ in $(\widehat{M},\mathcal{L})$ is simply denoted by
\begin{equation}\label{alg_count}
\mathcal{M}_S(\gamma;(y_{\zeta_1},\cdots,y_{\zeta_k})):=
\widetilde{\mathcal {M}}_S(\gamma;(y_{\zeta_1},\cdots,y_{\zeta_k}); \widehat{M},\mathcal{L})/\mathrm{Aut}(S).
\end{equation}
And denote simply by
$$\mathcal{M}(\gamma;(y_{\zeta_1},\cdots,y_{\zeta_k}))=
\bigcup_{S\in\mathcal{R}^1_k}\mathcal{M}_S(\gamma;(y_{\zeta_1},\cdots,y_{\zeta_k}))$$
the total moduli space of $(1,k)$-punctured pseudo-holomorphic
disks in $(\widehat{M},\mathcal{L})$.

As we only consider exact case, in this paper we ignore all transversality problems
and assume that regularity is satisfied for all involved moduli spaces. Alternatively,
we need the assumption that the adjusted almost complex structure $J\in\mathcal
{J}(\lambda,\widehat{\omega})$ is {\it
regular} which means that the linearized operator $D_u$ is
surjective and  the
moduli space of $(1,k)$-punctured pseudo-holomorphic maps from $S$ ($k\geq 1$) has expected dimension
\begin{eqnarray}\label{dim-moduli}
\dim\mathcal{M}_S(\gamma;(y_{\zeta_1},\cdots,y_{\zeta_k}))
 &=&\mu_{CZ}(\gamma)
-\deg(y_{\zeta_1})-\cdots-\deg(y_{\zeta_k}).
\end{eqnarray}

\subsection{Compactification}

Denote by ${S^\circ=S-{\rm Im}(\epsilon_\eta)}$. For $u\in \widetilde{\mathcal {M}}
(\gamma;(y_{\zeta_1},\cdots,y_{\zeta_k}); \widehat{M},\mathcal{L})$, we define its energy by
$$E(u)=E_H(u\circ\epsilon_\eta)+E(u|_{S^\circ}),$$
where $E_H(u\circ\epsilon_\eta)$ is the Hofer energy (see \cite{BEHWZ}) of the
pseudo-holomorphic negative-cylinder $F=u\circ\epsilon_\eta$ asymptotic to a closed Reeb orbit,
and $E(u|_{S^\circ}):=\int_{S^\circ}\frac{1}{2}|Du-Y|^2$ is the analogue of the energy
defined in \cite[(8.12)]{Seidel}.

Since $\widehat{M}$ is an  almost complex manifold with symmetric cylindrical
ends adjusted to the symplectic form $\widehat{\omega}$, by applying
Proposition 6.3 of \cite{BEHWZ} to $u\circ\epsilon_\eta$ and the usual argument of
Floer homology to $u|_{S^\circ}$, the energies $E(u)$ are uniformly bounded for all
$u\in \widetilde{\mathcal {M}}(\gamma;(y_{\zeta_1},\cdots,y_{\zeta_k}))$. Thus,  by Theorem 10.2 of
\cite{BEHWZ} and the usual Floer's compactifying method of involving broken curves,
one can obtain the compactification moduli space
$\overline{\mathcal {M}}(\gamma;(y_{\zeta_1},\cdots,y_{\zeta_k}))$.

\begin{theorem}[Compactness Theorem]\label{comp-thm}
Let $M$ be an exact symplectic manifold with contact type boundary $W$
such that $c_1(M)=0$,
$\widehat{M}$ be its symplectic completion, $J\in\mathscr J$, and $\mathcal L=\{L_C\}$ be a
collection of closed admissible Lagrangian submanifolds in $M$ which do not intersect
with $W$. Let $\gamma\in\mathcal P_\lambda$ and let
$(y_1,\cdots,y_k)$ be a collection of
Hamiltonian chords decided by the chosen Floer data and perturbation data.
Then for generic choices of consistent universal analytic data {\bf D}$_{rel}$, 
the moduli space $\overline{\mathcal {M}}(\gamma;(y_1,\cdots,y_k))$
is a smooth compact stratified manifold of dimension
\begin{eqnarray}\label{dim_of_moduli}
\dim\mathcal{M}(\gamma;(y_{1},\cdots,y_{k})) &=& |\gamma|-n+k+2-\deg(y_1)-\cdots-\deg(y_k),
\end{eqnarray}
whose codimension one strata consist of moduli spaces of the following form
\begin{eqnarray}\label{comp-boundary}
&&{\mathcal {M}}(\gamma;\gamma',\gamma_1,\cdots,\gamma_n)/\mathbb R\times{\mathcal{M}}(\gamma';y_1,\cdots,y_k)
\times
\mathcal{M}(\gamma_1;\emptyset)\times\cdots\times\mathcal{M}(\gamma_n;\emptyset)
\nonumber\\
&\bigcup&
{\mathcal{M}}(\gamma; y_1,\cdots,y_{i-1},x,y_j\cdots,y_{k})\times{\mathcal{M}}( x,y_{i+1}\cdots,y_{j-1}),
\end{eqnarray}
where $\gamma'$ runs over all possible closed Reeb orbits and
$x$ is any element in, say, $\Hom(\tilde L_i,\tilde L_j)$.
\end{theorem}

\begin{proof}
The first type of boundary strata
comes from the neck-stretching compactification process (\cite{BEHWZ})
and the second type of boundary strata
arises from the disk-bubbling-off compactification in Fukaya category (\cite[\S9]{Seidel}).
\end{proof}

\subsection{Compactification for case of several punctures}

For the case of several interior punctures, the
boundary of the compactified moduli space becomes
slightly more complicated; namely, when compactifying
the moduli spaces, the neck-stretching of pseudo-holomorphic
curves will produce not only pseudo-holomorphic planes
in $\widehat M$ (or say pseudo-holomorphic spheres with one puncture
at $+\infty$), but also pseudo-holomorphic spheres with several punctures
at $+\infty$ (compare with formula \ref{def_bracket} for the Lie bracket
of the linearized contact chain complex). We state the theorem below, and
leave the proof to the interested reader:

\begin{theorem}\label{comp-thm-2punctures}
Assume the conditions in Theorem \ref{comp-thm}. Let
$\gamma_1^+,\gamma_2^+,\cdots,\gamma_r^+\in\mathcal P_\lambda$ and let $(y_1,\cdots,y_k)$ be a
collection of Hamiltonian chords decided by the chosen Floer data and perturbation data.
Then  for generic choices of consistent universal analytic data {\bf D}$_{rel}$, the
moduli space $\overline{\mathcal {M}}(\gamma_1^+,\cdots,\gamma_r^+;(y_1,\cdots,y_k))$
is a smooth compact stratified manifold of dimension
\begin{equation}\label{dim_of_moduli-2}
\dim\mathcal{M}(\gamma_1^+,\cdots,\gamma_r^+;(y_{1},\cdots,y_{k}))
= k+2r-3+\mu(\gamma_1^+)+\cdots+\mu(\gamma_r^+)-\deg(y_1)-\cdots-\deg(y_k),
\end{equation}
whose codimension one strata are the union of moduli spaces of broken (or say, neck-stretching)
curves that appear in SFT and moduli spaces of disk-bubbling-off curves that appear in the Fukaya
category; more precisely, they consists of 
\begin{enumerate}
\item the products of moduli spaces $\mathcal M(\gamma_1^+,\cdots,
\hat\gamma_{p_1}^+,\cdots,\hat\gamma_{p_s}^+,\cdots,\gamma_r^+;
\gamma^-,\gamma_1^-,\cdots,\gamma_n^-)/\mathbb R$, where ``$\hat\quad$'' 
means the corresponding item is omitted, together with moduli spaces of vertical cylinders over each $\gamma_{p_j}^+$,
and $\mathcal M(\gamma^-;(y_{1},\cdots,y_{k}))$ 
together with 
$\mathcal M(\gamma_{q_1}^+,\cdots,\gamma_{q_t}^+, \gamma_j^-;\emptyset)$, 
for some $\{\gamma_{q_1}^+,\cdots,\gamma_{q_t}^+\}$ a subset of\linebreak $\{\gamma_{p_1}^+,\cdots,
\gamma_{p_s}^+\}$;

\item the products of moduli spaces 
$\mathcal M(\gamma_1^+,\cdots,\hat\gamma_{p_1}^+,\cdots,\hat\gamma_{p_s}^+,
\cdots,\gamma_r^+;( y_1,\cdots,y_{i-1},x,y_j,\cdots,y_{k}))$ and 
$\mathcal M(\gamma_{p_1}^+,\cdots,\gamma_{p_s}^+;( x,y_{i},\cdots,y_{j-1}))$, for some
$x\in\Hom(\tilde L_i,\tilde L_j)$.
\end{enumerate}
\end{theorem}

\subsection{Orientation}

If $\mbox{Char}(\mathbb K)=2$, the regularity and compactness of those involved moduli spaces will be
enough to define the relating maps. For an arbitrary field $\mathbb K$,
we need to consider the orientation problem of the moduli spaces. 
The algorithm to orient the 
moduli spaces is standard nowadays, and the interested reader may refer
to, for example, \cite{BM,FOOO,H,Seidel,So}
for discussions in various situations with respect to punctured or/and bordered
Riemann surfaces.

Basically, one can show that the linearization of the inhomogeneous pseudo-holomorphic map equation
gives rise to a Fredholm operator between two corresponding Banach spaces, and then apply the construction of 
coherent orientations for Fredholm operators to moduli spaces
of holomorphic maps in $\widehat M$ relating closed Reeb orbits in $W$ and Hamiltonian chords in $M$.

The paper of Bourgeois-Mohnke \cite{BM} is of particular interest to us,
since we may simply replace in their paper the orientation bundle of one of the closed Reeb orbits with
(the product of) the orientation of Hamiltonian chords, and all their arguments can be applied to our case.

As a concrete example, suppose $u:S\to\widehat M$ is a solution to equation (\ref{puncture-ipm}).
Then the linearized operator $D_u$ corresponding to (\ref{puncture-ipm}) 
is a Fredholm operator, whose determinant line bundle $\det D_u$ is
isomorphic to $\det\mathcal (R_k^1)_u\otimes\bigotimes_i o_{y_i}\otimes o_{\gamma}^{-}$.
On the other hand, both in SFT and in Fukaya category, we have a gluing map (\cite[\S3]{BM})
\begin{equation}\label{glue-diff}
\mathcal{G}: \mathcal{M}(\gamma^+;\gamma^-; W\times\R)\times \mathcal {M}
(\gamma^-;(y_1,\cdots,y_k); \widehat{M},\mathcal{L})
 \to \mathcal{M}(\gamma^+;(y_1,\cdots,y_k); \widehat{M},\mathcal{L}),
\end{equation}
and a gluing map (\cite[\S12]{Seidel})
\begin{equation}\label{glue-diff-bound}
\mathcal {F}: \mathcal{M}(\gamma;(y_1,\cdots,y_l,y); \widehat{M},
\mathcal{L}) \times \mathcal{M}((y,y_{l+1},\cdots,y_k); \widehat{M},\mathcal{L}')
\to \mathcal{M}(\gamma;(y_1,\cdots, y_k); \widehat{M},\mathcal{L}\cup\mathcal{L}').
\end{equation}
Thus to obtain an orientation for $\mathcal M(\gamma^+;(y_1,\cdots,y_k); \widehat M,\mathcal L)$,
one choose step by step from the orientations of the moduli spaces
that arise in SFT and Fukaya category respectively such that
the gluing maps (\ref{glue-diff}) and (\ref{glue-diff-bound}) are orientation preserving.

In general, for the case of several punctures, we state the following theorem without proof:

\begin{theorem}\label{Orient-glue}
Under the conditions of Theorem \ref{comp-thm},
the determinant
line bundles over the moduli spaces 
\begin{eqnarray*}
\mathcal{M}(\gamma_1^+,\cdots,\gamma_r^+;(y_1,\cdots,y_k);\widehat M,\mathcal L\cup\mathcal L'),&&
\mathcal {M}(\gamma_1^+,\cdots,\gamma_r^+;\gamma^-; W\times\R),\\\mathcal {M}
(\gamma^-;(y_1,\cdots,y_l, y); \widehat{M},\mathcal{L}) ,&&
\mathcal{M}((y,y_{l+1},\cdots,y_k); \widehat{M},\mathcal{L}')\end{eqnarray*}
are orientable in such a way that the gluing 
maps $\mathcal{G}$ and $\mathcal {F}$ preserve the orientations up to sign.
Moreover, if we switch $\gamma^+_i$ with $\gamma^+_{i+1}$, then
the orientation of $\mathcal M(\gamma^+_1,\cdots,\gamma^+_r;(y_1,\cdots,y_k);\widehat M,\mathcal L\cup\mathcal L')$ 
changes its sign by $(-1)^{|\gamma^+_i||\gamma^+_{i+1}|}$; and
if we cyclically permute $y_1,\cdots,y_k$, then its orientation changes its sign by
$(-1)^{|y_k|(|y_1|+\cdots+|y_{k-1}|)}$.
\end{theorem}

\section{Homomorphism from linearized contact homology to cyclic cohomology}\label{sect_rel}

\subsection{The chain homomorphism}
We are now ready to give a chain homomorphism from the linearized chain complex of $M$ to
the cyclic cochain complex of the Fukaya category of $M$.

\begin{theorem}\label{thm_liebimap}
Let $M$ be an exact simply connected symplectic manifold with contact type boundary such that $c_1(M)=0$.
Let
$$
f:\Cont_*^{\mathrm{lin}}(M)\to \Cycl^*(\Fuk(M))
$$
be the homomorphism
such that for each $\gamma\in\Cont_*^{\mathrm{lin}}(M)$,
the value of $f(\gamma)$ at $(\ov y_1,\ov y_2,\cdots,\ov y_m)$
is
$$f(\gamma)(\ov y_1,\ov y_2,\cdots,\ov y_m):=
\#\mathcal M(\gamma;(y_1,y_2,\cdots, y_m)),$$
where $\mathcal M(-)$ is given in (\ref{alg_count}).
Then $f$ thus defined is a chain homomorphism.
\end{theorem}

\begin{proof}
For $\gamma\in\Cont_*^{\mathrm{lin}}(W)$, and $(\ov y_1,\ov y_2,\cdots,\ov y_m)\in\ov{\Hom}(\tilde L_1,\tilde L_2)\otimes
\ov{\Hom}(\tilde L_2,\tilde L_3)\otimes\cdots\otimes\ov{\Hom}(\tilde L_m,\tilde L_1)$,
suppose the moduli space $\mathcal M(\gamma;(y_1, y_2,\cdots, y_m))$ is 1-dimensional, then
by the compactness theorem (Theorem \ref{comp-thm}) its compactified boundary is composed
of the following two types of broken pseudo-holomorphic curves:
\begin{itemize}
\item The first type consists of broken curves
which are stretched from the pseudo-holomorphic disks in $W\times\mathbb R^+$.
During the stretching, the upper cylinder may generate some pseudo-holomorphic planes
in $M$;

\item The second type consists of broken curves
which includes the bubbling-off disks of the original pseudo-holomorphic disks.
\end{itemize}
To show $f(d\gamma)=b(f(\gamma))$, consider their values at
$(\ov y_1,\ov y_2,\cdots,\ov y_m)$:
$f(d\gamma)(\ov y_1,\ov y_2,\cdots,\ov y_m)$ is exactly the number of broken pseudo-holomorphic disks of the first type
and $b(f(\gamma))(\ov y_1,\ov y_2,\cdots,\ov y_m)$, which equals $f(\gamma)(b(\ov y_1,\ov y_2,\cdots,\ov y_m))$,
is exactly the number of broken pseudo-holomorphic
disks of the second type.
That means, $f$ is a chain map for at least the case
$\mbox{Char}(\mathbb K)=2$.

For general $\mathbb K$, this still holds since as we have commented before, the orientations of
the moduli space of the pseudo-holomorphic disks in Fukaya category
and the moduli space of the pseudo-holomorphic punctured spheres in SFT both
follows the Koszul sign rule, and therefore the algebraic counting of the elements
in two types of boundary strata are equal.
\end{proof}

\subsection{Homomorphism of Lie$_\infty$ algebras}

We have known from previous sections
that the linearized contact chain complex forms a Lie$_\infty$ algebra,
and the cyclic cochain complex of the Fukaya category is a Lie algebra, which
is a special class of Lie$_\infty$ algebras.
In this subsection we are going to show that the chain map defined in Theorem
\ref{thm_liebimap} is in fact a Lie$_\infty$ algebra homomorphism.

We first introduce the definition of the homomorphism of Lie$_\infty$ algebras.

\begin{definition}[Lie$_\infty$ Homomorphism]
Suppose $V$ and $W$ are two Lie$_\infty$ algebras.
A Lie$_\infty$ homomorphism from $V$ to $W$ is a differential graded coalgebra map
$$F:\mbox{$\bigwedge^\bullet V$}\longrightarrow\mbox{$\bigwedge^\bullet W$}.$$
More precisely, a Lie$_\infty$ homomorphism
from $V$ to $W$ consists of a sequence of linear operators
$$F_k:\mbox{$\bigwedge^k V$}\longrightarrow W,\quad k=1,2,\cdots$$
such that for all $a_1a_2\cdots a_k\in\bigwedge^k V$,
\begin{multline}\label{Lieinfty_map}
\sum_i\sum_{k_1+\cdots+k_i=k}\sum_{\sigma}
\delta_i(F_{k_1}(a_{\sigma(1)}\cdots a_{\sigma(k_1)})\cdots F_{k_i}
(a_{\sigma(k_1+\cdots+k_{i-1}+1)}\cdots a_{\sigma(k)}))\\
=\sum_{p+q=k+1}\sum_{\mu} F_p(\delta_q(a_{\mu(1)}\cdots a_{\mu(q)})a_{\mu(q+1)}\cdots a_{\mu(k)}),
\end{multline}
where $\sigma$ runs over all $(k_1,k_2,\cdots,k_i)$-unshuffles of $k$, and
$\mu$ runs over all $(q,p)$-unshuffles of $k$, respectively.
\end{definition}

A Lie$_\infty$ homomorphism from $V$ to $W$ induces a Lie algebra homomorphism
from $\HH_*(V,\delta_1)$ to $\HH_*(W,\delta_1)$. We have:

\begin{theorem}\label{thm_liebimap_comp}
Let $f$ be the map in Theorem~\ref{thm_liebimap}. Then $f$ can be completed to
be a Lie$_\infty$ homomorphism, {\it i.e.} there is a sequence of operators $f_1=f,
f_2,\cdots,$ satisfying $(\ref{Lieinfty_map})$, where
$$f_k:\mbox{$\bigwedge^k$}\Cont_*(M)\longrightarrow \Cycl^*(\Fuk(M)),\quad k=1,2,\cdots.$$ 
In particular, $f$ induces a Lie bialgebra map
from the linearized contact homology of $M$ to the cyclic cohomology of the Fukaya
category of $M$.
\end{theorem}

\begin{proof}
Let $f_1$ be the map $f$ in Theorem~\ref{thm_liebimap}.
For $k>1$, define
$$f_k: \mbox{$\bigwedge^k$}\Cont_*(M)\to\Cycl^*(\Fuk(M))$$ as follows:
$$
f_k(\gamma_1\gamma_2\cdots\gamma_k)(\ov y_1,\ov y_2,\cdots,\ov y_m)=
\#\mathcal M(\gamma_1,\gamma_2,\cdots,\gamma_k;(y_1, y_2,\cdots, y_m)).
$$
We show that $f_1,f_2,\cdots$ thus defined is a homomorphism of Lie$_\infty$ algebras.

In fact, consider the moduli space
$\mathcal M(\gamma_1,\gamma_2,\cdots,\gamma_k;(y_1,y_2,\cdots,y_m))$.
If it is of dimension one, then the boundary of the compactification of $\mathcal M$
consists of two types of pseudo-holomorphic curves, with ``neck stretching'' and ``disk bubbling-off'',
respectively:
\begin{itemize}
\item[(I)] the stretching occurs in $M\times \mathbb R^+$, and all such possibilities correspond to
the right hand side of equation (\ref{Lieinfty_map});

\item[(II)] the bubbling-off occurs in $W$, which consists of all possibilities of disk bubbling-off,
(the case that one of the saddle point in the Riemann surface is pushed down to zero is included),
and all all such possibilities correspond to the left hand side of equation (\ref{Lieinfty_map}).
\end{itemize}

As a concrete example, we consider the moduli space of pseudo-holomorphic disks with two punctures, as in
the following picture:
\begin{center}
\includegraphics[height=3.6cm]{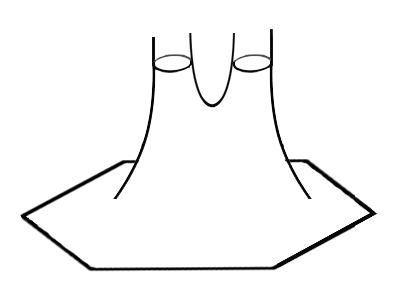}
\end{center}
where $\gamma_1,\gamma_2\in\Cont_*^{\mathrm{lin}}(W)$ are two punctures and
the boundary of the disks lies in $\tilde L_1\cup \tilde L_2\cup\cdots\cup \tilde L_m$ in cyclic order.
Suppose the moduli space is 1-dimensional, then its boundary points
consist of four types of broken pseudo-holomorphic
curves:
\begin{enumerate}
\item The stretching of the pants occurs at one of the upper sleeves but not the other; in other words, one
of them is stretched to be infinitely long; during the stretching, this sleeve may generate some sphere bubbles in
$M$, which can be capped off, however, during the stretching, the other sleeve remain unchanged.

\item The stretching of the pants occurs above the bottom sleeve;
it consists of two sub cases, one is that during the stretching, the upper half of the pants may
generate some sphere bubbles in $M$ which will be capped off; and the other is that,
during the stretching, one of the sleeves remains unchanged, and it will be capped off in $M$ with the bubble
generated by the other sleeve;

\item The bottom sleeve splits into two pieces, with each piece has a puncture at $+\infty$ in it,
{\it i.e.} the saddle point in the pants is pushed down into $M$ which splits the pants along a intersection point
of two Lagrangian submanifolds;

\item The bottom sleeve has a disk bubbling-off.
\end{enumerate}
These four cases correspond to the following operations respectively:
(1) $f_2(d\gamma_1,\gamma_2)$ and $f_2(\gamma_1,d\gamma_2)$; 
(2) $f_1[\gamma_1,\gamma_2]$; (3) $[f_1(\gamma_1),f_1(\gamma_2)]$; 
and (4) $b\circ f_2(\gamma_1,\gamma_2)$.
It follows from $(1)+(2)=(3)+(4)$ (by Compactness Theorem \ref{comp-thm-2punctures})
that $$[f_1(\gamma_1),f_1(\gamma_2)]-f_1[\gamma_1,\gamma_2]=f_2(d\gamma_1,\gamma_2)+
f_2(\gamma_1, d\gamma_2)-b\circ f_2(\gamma_1,\gamma_2).$$
This is exactly equation (\ref{Lieinfty_map}) for $k=2$, {\it i.e.} $f_1$ maps Lie bracket to Lie bracket up to homotopy.
For the general case, the argument is similar.
\end{proof}

\subsection{Lie bialgebras and more}

As we have said before, Cieliebak-Latschev have shown that the linearized contact complex
endows the structure of a BV$_\infty$ algebra. We remark that the 
cyclic cochain complex of the Fukaya category naturally endows a BV$_\infty$ algebra structure
in the sense of Cieliebak-Latschev, and the homomorphism
in Theorem \ref{thm_liebimap_comp} is in fact a BV$_\infty$ algebra homomorphism.
For example, to show that $f$ maps the cobracket to the cobracket up to homotopy, we consider
pseudo-holomorphic cylinders with two components of the boundary each lying in one set
of cyclic chain complex in $\mathcal Fuk(M)$ and with
one puncture at $+\infty$ (or say, the pseudo-holomorphic
pants with two sleeves lying in two sets of Lagrangian submanifolds and one sleeve going to $+\infty$), with
the similar argument as above, one sees that $f$ is a Lie coalgebra map up to homotopy.

\section{Example of cotangent bundles}\label{sect_cot}

In this last section we briefly discuss the cyclic homology of Fukaya category
in cotangent bundles. In general, the Fukaya category is very difficult
to compute. However, in the case of cotangent bundles, the Fukaya
category is strikingly simple, due to the following theorem:

\begin{theorem}[Fukaya-Seidel-Smith and Nadler]\label{thm_FSS}
Let $N$ be a simply connected, compact spin manifold and $T^*N$ its cotangent bundle.
Then the Fukaya category of $T^*N$ is derived Morita equivalent to the zero section $N$.
\end{theorem}

This theorem is independently and simultaneously proved by Fukaya-Seidel-Smith in \cite[Theorem 1]{FSS}
and Nadler in \cite[Theorem 1.3.1]{Nadler}. A further discussion of this result is given in Fukay-Seidel-Smith \cite{FSS2}.

In fact, what Fukaya et. al. proved is even stronger than the above theorem, where an explicit homotopy equivalence
of two categories is given. Namely,
for two admissible Lagrangian submanifolds $\tilde L_1, \tilde L_2$, there is a $c\in\Hom^0(\tilde L_1, \tilde L_2)$, such
that
$$\xymatrix{
\Hom(\tilde L_1,\tilde L_1)\ar[r]^{c\cup\mbox{-}}&\Hom(\tilde L_1,\tilde L_2),
}$$
is a quasi-isomorphism,
where $c\cup\mbox{-}$ means the composition with $c$ (see the last paragraph of \cite[\S1]{FSS}).

A general theory in homological algebra says that
if two DG categories are derived Morita equivalent, then their cyclic homology groups
are isomorphic ({\it c.f.} To\"en \cite[\S5.2]{Toen2}). Indeed, To\"en proved the isomorphism
of Hochschild (co)homology groups; the isomorphism of cyclic homology groups
can then be obtained by comparing the associated
Connes' long exact sequences, or by showing that the cyclic homolology is a derived functor which 
is invariant under derived equivalences. This property holds for A$_\infty$ categories as well, since any A$_\infty$
category is homotopy equivalent to a DG category ({\it c.f.} Fukaya \cite[Corollary 9.4]{Fukaya}).
As a corollary, we have the following theorem:

\begin{theorem}\label{CorFSS}
Let $N$ be a
simply-connected, compact spin manifold
and $T^*N$ the cotangent bundle of $N$.
Then the cyclic homology of $\mathcal Fuk(T^*N)$ is
isomorphic to the cyclic homology of the Floer cochain complex
$CF^*(\tilde N)$.
\end{theorem}

On the other hand, the Floer cochain complex
$CF^*(\tilde N)$, and even its cyclic homology, is known for symplectic geometers
by the following two theorems:

\begin{theorem}[The PSS Isomorphism]\label{PSS}
Let $N$ be a simply-connected manifold as before.
Then the Floer cochain complex $CF^*(\tilde N)$ of $N$ is quasi-isomorphic to
its de Rham cochain complex $\Omega^\bullet(N)$.
\end{theorem}

\begin{proof}This is the chain level statement of the famous Piunikhin-Salamon-Schwarz (PSS
for short) isomorphism,
first appeared in \cite{PSS}.
A proof in the most general case (without assuming the result of Fukaya-Seidel-Smith and Nadler cited above),
can be found in Abouzaid \cite[Theorem 1.1]{Ab092}.
\end{proof}

\begin{theorem}[K.-T. Chen and Jones]\label{freeloopiso}
Let $N$ be a simply connected manifold and $LN$ its free loop space.
Denote by $\Omega^\bullet(N)$ the de Rham cochain complex of $N$. Then
the cyclic homology of $\Omega^\bullet(N)$ is isomorphic to the equivariant cohomology $\HH^*_{S^1}(LN)$.
\end{theorem}

\begin{proof}
See K.-T. Chen \cite[Theorem 4.3.1]{KTChen} and Jones \cite[Theorem A]{Jones}.
\end{proof}

Combining Theorems \ref{CorFSS}, \ref{PSS} and \ref{freeloopiso} yields the following:

\begin{theorem}Let $N$ be a
simply-connected spin manifold
and $T^*N$ the cotangent bundle of $N$.
Then the cyclic cohomology of $\mathcal Fuk(T^*N)$ is isomorphic to the
equivariant homology $\HH_*^{S^1}(LN)$ of the free loop space $LN$.
\end{theorem}

\subsection{Relations to the result of Cieliebak-Latschev}
In the article \cite{CL} Cieliebak-Latschev have shown
(see \cite[Theorem C]{CL}) that
the linearized contact homology of $T^*N$ is isomorphic, as involutive Lie bialgebras, to the
relative $S^1$-equivariant homology
$\HH_*^{S^1}(LN,N)$, where $N$ is identified with the set of constant loops.
The Lie bialgebra of the later is obtained by Chas-Sullivan in string topology (\cite{CS02}).
The map of Cieliebak-Latschev is also given by considering
the pseudo-holomorphic cylinders with one boundary approaching a closed
Reeb orbit and the other lying on a Lagrangian submanifold (zero section in this case).
There is a long exact sequence (of excision)
\begin{equation}\label{MV}
\xymatrix{
\cdots\ar[r]&\HH_*^{S^1}(N)\ar[r]^{i_*}&\HH_*^{S^1}(LN)\ar[r]^{j_*}&\HH_*^{S^1}(LN,N)\ar[r]&\HH_{*-1}^{S^1}(N)\ar[r]&\cdots
}
\end{equation}
for the pair $(LN,N)$.
This suggests the following diagram
$$
\xymatrixcolsep{3.5pc}
\xymatrixrowsep{2.5pc}
\xymatrix{
&&\mathrm{CH}^{\mathrm{lin}}_*(T^*N)\ar[lld]_-{\mbox{\scriptsize
Theorem~\ref{mapfromcont}}}^-{f}\ar[ld]^-{f|N}\ar[d]\ar[rd]_-{\cong}^-{\mbox{\scriptsize Cieliebak-Latschev}}&\\
\mathrm{HC}^*(\Fuk(T^*N))\ar[r]^-{\cong}_-{\mbox{\scriptsize Morita equiv.}}
&\mathrm{HC}^*(CF^*(\tilde N))\ar[r]^-{\cong}_-{\mbox{\scriptsize Chen-Jones}}
&\mathrm{H}_*^{S^1}(LN)\ar[r]^{j_*}_-{(\ref{MV})}&\mathrm{H}_*^{S^1}(LN,N).
}
$$
However, the diagram may in general not be commutative between the left and right, since otherwise the long
exact sequence in (\ref{MV}) will be splitting, which is not true in general (pointed to us by Pomerleano).

\end{document}